\documentclass[final,authoryear,1p,times]{elsarticle}
\usepackage{graphicx,subfig}
\usepackage{amsmath,amssymb}
\usepackage{url}
\usepackage{algorithm}
\usepackage{algorithmic}

\begin{document}

\begin{frontmatter}

\title{Adaptive discontinuous Galerkin methods for non-linear diffusion-convection-reaction equations}

\author[mu]{Murat Uzunca\corref{cor}}
\ead{uzunca@metu.edu.tr}

\author[bk]{B\"{u}lent Karas\"{o}zen}
\ead{bulent@metu.edu.tr}

\author[mm]{Murat Manguo\u{g}lu}
\ead{manguoglu@ceng.metu.edu.tr}

\cortext[cor]{Corresponding author. Tel.: +90 312 2105353.}

\address[mu]{Department of Mathematics, Middle East Technical University, 06800 Ankara, Turkey}
\address[bk]{Department of Mathematics and Institute of Applied
  Mathematics, Middle East Technical University, 06800 Ankara, Turkey}
\address[mm]{Department of Computer Engineering and Institute of Applied
  Mathematics, Middle East Technical University, 06800 Ankara, Turkey}

\begin{abstract}
In this work, we apply the adaptive discontinuous Galerkin (DGAFEM) method to the convection dominated non-linear, quasi-stationary diffusion-convection-reaction equations.
 We propose  an efficient preconditioner using a matrix reordering scheme to solve the sparse linear systems iteratively arising from the discretized non-linear equations. Numerical examples demonstrate effectiveness of the DGAFEM to damp the spurious oscillations and resolve well the sharp layers occurring in convection dominated non-linear equations.
\end{abstract}


\end{frontmatter}

\section{Introduction}
Many engineering problems such as chemical reaction processes, heat conduction,
 nuclear reactors, population dynamics etc. are governed by coupled convection-diffusion-reaction partial
differential equations (PDEs) with non-linear source or sink terms. It is a significant challenge to solve such PDEs numerically when they are convection/reaction-dominated. As a model problem, we consider the coupled quasi-stationary  equations arising from the time discretization of time-dependent non-linear diffusion-convection-reaction equations \citep{bause12ash}
\begin{subequations}\label{1}
\begin{align}
  \alpha u_i-\epsilon_i\Delta u_i + {\bf b}_i\cdot\nabla u_i + r_i({\bf u}) &= f_i \quad \, \text{ in } \; \Omega_i, \\
   u_i &= g^D_i  \quad \text{on } \; \Gamma^D_i, \\
	\epsilon_i\nabla u_i\cdot {\bf n} &= g^N_i  \quad \text{on } \; \Gamma^N_i, \quad  i = 1,\ldots, m
\end{align}
\end{subequations}
with $\Omega_i$ are bounded, open, convex domains in $\mathbb{R}^2$ with boundaries $\partial\Omega_i =\Gamma^D_i\cup\Gamma^N_i$, $\Gamma^D_i\cap\Gamma^N_i=\emptyset$, $0<\epsilon_i\ll 1 $ are the diffusivity constants, $f_i\in L^2(\Omega)$ are the source functions, ${\bf b}_i\in\left(W^{1,\infty}(\Omega)\right)^2$ are the velocity fields, $g^D_i\in H^{3/2}(\Gamma^D_i )$ are the Dirichlet boundary conditions, $g^N_i\in H^{1/2}(\Gamma^N_i )$ are the Neumann boundary conditions, and ${\bf u}(x)=(u_1, \ldots , u_m)^T$ and ${\bf n}$ denote the vector of unknowns and normal vector to the boundary, respectively. The coefficients of the linear reaction terms, $\alpha >0$, stand for the temporal discretization, corresponding to $1/\Delta t$, where $\Delta t$ is the discrete time-step. Moreover, we assume that the non-linear reaction terms are bounded, locally Lipschitz continuous and monotone, i.e. satisfy for any $s, s_1, s_2\ge 0$, $s,s_1, s_2 \in \mathbb{R}$ the following conditions  \citep{bause12ash}
\begin{subequations}
\begin{align}
|r_i(s)| &\leq C , \quad C>0\\
\| r_i(s_1)-r_i(s_2)\|_{L^2(\Omega)} &\leq L\| s_1-s_2\|_{L^2(\Omega)} , \quad L>0\\
 r_i\in C^1(\mathbb{R}_0^+), \quad r_i(0) =0, &\quad r_i'(s)\ge 0.
\end{align}
\end{subequations}
The non-linear reaction terms $r_i({\bf u})$ occur in chemical engineering usually in the form of products and rational functions of concentrations, or exponential functions of the temperature, expressed by the Arrhenius law.
Such models describe chemical processes and they are strongly coupled as an inaccuracy in one unknown affects all the others.  Hence, efficient numerical approximation of these systems is needed. For the convection/reaction-dominated problems, the standard Galerkin finite element methods are known to produce spurious oscillations, especially in the presence of sharp fronts in the solution, on boundary and interior layers.

In last two decades several stabilization and shock/discontinuity
capturing techniques were developed for linear and non-linear stationary and time dependent
problems.  For linear convection dominated problems, the streamline upwind Petrov-Galerkin(SUPG) method is capable to stabilize the unphysical oscillations \citep{bause10sfe, bause12ash, bause13hof}. Nevertheless, in non-linear convection dominated problems, spurious oscillations are still present in crosswind direction. Therefore, SUPG is used with the anisotropic shock capturing technique as SUPG-SC for reactive transport problems \citep{bause10sfe, bause12ash, bause13hof}. It was shown that SUPG-SC is capable of reducing the unphysical oscillations in cross-wind direction. The parameters of the SUPG and SUPG-SC should be designed carefully for the efficient solution of the discretized equations.

In contrast to the  standard Galerkin conforming finite element methods, discontinuous Galerkin finite element(DGFEM) methods produce stable discretizations without the need for stabilization strategies, and damp the unphysical oscillations for linear convection dominated problems. In \citep{yucel13dgf}, several non-linear convection dominated problems of type (\ref{1}) are solved with DG-SC, discontinuous Galerkin method with the shock-capturing technique in \citep{persson06ssc}. For an accurate solution of non-linear convection dominated problems, higher order finite elements are used because they are less diffusive and avoid artificial mixing of chemical species under discretization, for SUPG-SC and DG-SC, respectively. The main advantages of DGFEM are the flexibility in handling non-matching grids and in designing hp-refinement strategies \citep{houston02dhp}, which allow easily adaptive grid refinement. In this paper we extend the adaptive discontinuous Galerkin method in  \citep{schotzau09rae} to the convection dominated non-linear problems of type (1).  We show on several examples the effectiveness and accuracy of  DGAFEM  capturing boundary and internal layers very sharply and without significant oscillations.  An important
drawback is that the resulting linear systems are  more dense than the continuous finite elements and ill-conditioned. The condition
number grows rapidly with the number of elements and with the penalty parameter. Therefore, efficient solution strategies  such as preconditioning are required to solve the linear systems. While more robust compared to iterative solvers, direct solvers are
usually more memory and time consuming due to fill-in. Furthermore, they
are known to be less scalable on parallel architectures. Therefore, in this
paper we use an iterative method which is robust and efficient.

Because the DG method requires more degrees of freedom than the standard Galerkin method, in \citep{cangiani13scd} linear diffusion-convection-reaction equations are discretized near the boundary and inter layers by the interior penalty DG method, away from the layers by continuous Galerkin method. It was shown that combination of both methods will not affect the stability of the DG method. Another important class of non-linear diffusion-convection equations are those with the non-linear convective term, including the viscous Burger's equation. In the recent years, several effective interior penalty DG methods with efficient time integrators and with space and time adaptivity for this class of problems are developed, see for example
\citep{doljsi08iipg, doljsi05dg, doljsi13hpdgfem}. An important class of non-linear convection-diffusion equations are the pellet equations, which model the intra-particle mass and heat transport in porous catalyst pellets. A comprehensive review of  weighted residual methods, the orthogonal collocation, Galerkin, tau and least squares methods
is given in \citep{solsvik13ewr} for solving  the the linear and non-linear pellet equations, where the methods are compared with respect to convergence of the residuals and computational efficiency.

The rest of this work is organized as follows. In the next two sections, we give the DG discretization and describe the residual based adaptivity for non-linear diffusion-convection-reaction problems. In Section 4, we have compared the DGFEM and DGAFEM  with another class of stabilized methods, the GLSFEM (Galerkin least square finite element method) for a linear convection dominated problem. Section \ref{linear} deals with an efficient solution technique for solving the linear system arising from the DG discretization. In Section \ref{numerical}, we demonstrate the effectiveness and accuracy  of DGAFEM for handling the sharp layers arising in several examples with different type of non-linear reaction mechanisms. The paper ends with some conclusions.

\section{Symmetric discontinuous interior penalty Galerkin (SIPG) discretization}
In this Section, we describe the DG discretization of the model problem (1). We begin with the classical weak formulation of the scalar equation ($m=1$) of (\ref{1}): find $u\in U$ such that

\begin{equation} \label{2}
\int_{\Omega}(\epsilon\nabla u\cdot\nabla v+{\bf b}\cdot\nabla uv+\alpha uv)dx+\int_{\Omega}r(u)vdx = \int_{\Omega}fvdx +\int_{\Gamma^N}g^Nv ds\; , \quad \forall v\in V
\end{equation}
where the solution space $U$ and the test function space $V$ are given by

\begin{equation*}
U= \{ u\in H^1(\Omega) \, : \; u=g^D \text{ on } \Gamma^D\}, \quad
V= \{ v\in H^1(\Omega) \, : \; v=0 \text{ on } \Gamma^D\}.
\end{equation*}
Under the assumptions given in the previous section, the problem (\ref{2}) has a unique solution in $U$. The next step of the classical (continuous) FEM is to find an approximation to the problem (\ref{2}) using a conforming, finite-dimensional subspace $V_h\subset V$, which requires that the space $V_h$ contains functions of particular smoothness (e.g. when $V=H_0^1(\Omega)$, then we choose $V_h\subset \{ v\in C(\overline{\Omega}) : \; v=0 \text{ on } \partial\Omega\}$). On the other hands, discontinuous Galerkin methods make it easy to use the non-conforming spaces, in which case the functions in $V_h\not\subset V$ are allowed to be discontinuous on the inter-element boundaries.

In this article, the discretization of the problem (\ref{1}) is based on the symmetric discontinuous interior penalty Galerkin (SIPG) method, a type of discontinuous Galerkin methods, for the diffusion part \citep{arnold02uad,riviere08dgm} and the upwinding  for the convection part \citep{ayuso09dgm,houston02dhp}.

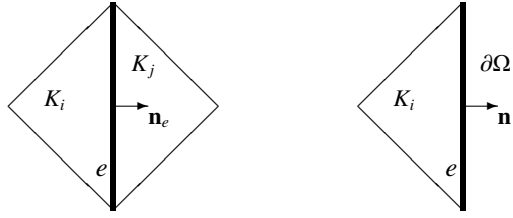
\begin{figure}
\centering
\setlength{\unitlength}{2.3mm}
\begin{picture}(35, 10)
\put(5,7){\line(1,1){6}}
\put(5,7){\line(1,-1){6}}
\put(17,7){\line(-1,1){6}}
\put(17,7){\line(-1,-1){6}}
\put(7,7){\small$K_{i}$ }
\put(12,9){\small$K_{j}$ }
\put(10,3){$e$ }
\put(11,7){\vector(1,0){2}}
\put(13,6){\small${\bf n}_{e}$}

\put(25,7){\line(1,1){6}}
\put(25,7){\line(1,-1){6}}
\put(27,7){\small$K_{i}$ }
\put(32,9){\small$\partial\Omega$ }
\put(30,3){$e$ }
\put(31,7){\vector(1,0){2}}
\put(33,6){\small${\bf n}$}

\linethickness{0.7mm} \put(11,1){\line(0,1){12}}
\linethickness{0.7mm} \put(31,1){\line(0,1){12}}
\end{picture}
\caption{Two adjacent elements sharing an edge (left); an element near to domain boundary (right)}
\label{jump}
\end{figure}

Let $\{\xi_h\}$ be a family of shape regular meshes with the elements (triangles) $K_i\in\xi_h$ satisfying $\overline{\Omega}=\cup \overline{K}$ and $K_i\cap K_j=\emptyset$ for $K_i$, $K_j$ $\in\xi_h$. Let us denote by $\Gamma^0$, $\Gamma^D$ and $\Gamma^N$ the set of interior, Dirichlet boundary and Neumann boundary edges, respectively, so that $\Gamma^0\cup\Gamma^D\cup\Gamma^N$ forms the skeleton of the mesh. For any $K\in\xi_h$, let $\mathbb{P}_k(K)$ be the set of all polynomials of degree at most $k$ on $K$. Then, set the finite dimensional solution and test function space by
$$
V_h=\left\{ v\in L^2(\Omega ) : v|_{K}\in\mathbb{P}_k(K) ,\; \forall K\in \xi_h \right\}\not\subset V.
$$
Note that the trial and test function spaces are the same because the boundary conditions in discontinuous Galerkin methods are imposed in a weak manner (see the SIPG construction below). Since the functions in $V_h$ may have  discontinuities along the inter-element boundaries, along an interior edge, there would be two different traces from the adjacent elements sharing that edge. In the light of this fact, let us first introduce some notations before starting the construction of SIPG formulation. Let $K_i$, $K_j\in\xi_h$ ($i<j$) be two adjacent elements sharing an interior edge $e=K_i\cap K_j\subset \Gamma_0$ (see Fig.\ref{jump}). Denote the trace of a scalar function $v$ from inside $K_i$ by $v_{i}$ and from inside $K_j$ by $v_{j}$. Then, set the jump and average values of $v$ on the edge $e$
$$
[v]= v_{i}{\bf n}_e- v_{j}{\bf n}_e , \quad \{ v\}=\frac{1}{2}(v_{i}+ v_{j}),
$$
where ${\bf n}_e$ is the unit normal to the edge $e$ oriented from $K_i$ to $K_j$. Similarly, we set the jump and average values of a vector valued function ${\bf q}$ on e
$$
[{\bf q}]= {\bf q}_{i}\cdot {\bf n}_e- {\bf q}_{j}\cdot {\bf n}_e , \quad \{ {\bf q}\}=\frac{1}{2}({\bf q}_{i}+ {\bf q}_{j}),
$$
Observe that $[v]$ is a vector for a scalar function $v$, while, $[{\bf q}]$ is scalar for a vector valued function ${\bf q}$. On the other hands, along any boundary edge $e=K_i\cap \partial\Omega$, we set
$$
[v]= v_{i}{\bf n} , \quad \{ v\}=v_{i}, \quad [{\bf q}]={\bf q}_{i}\cdot {\bf n}, \quad \{ {\bf q}\}={\bf q}_{i}
$$
where ${\bf n}$ is the unit outward normal to the boundary at $e$.

Now, we are ready to construct the SIPG discretization of the diffusion part of the problem. We start with the continuous (i.e. assume for the moment that the unknown solution $u$ is smooth enough) Poisson problem
\begin{subequations}
\begin{align}
  -\Delta u &= f \quad \;\; \text{ in } \; \Omega, \label{poisson} \\
   u &= g^D  \quad \text{ on } \; \Gamma^D \\
	\nabla u\cdot {\bf n} &= g^N  \quad \text{ on } \; \Gamma^N
\end{align}
\end{subequations}
We multiply the equation (\ref{poisson}) by a test function $v\in V_h$, we integrate over $\Omega$ and we split the integrals:
$$
-\sum \limits_{K \in {\xi}_{h}} \int_{K} \Delta uv dx = \sum \limits_{K \in {\xi}_{h}} \int_{K} fvdx
$$
Apply the divergence theorem on every element integral
$$
\sum \limits_{K \in {\xi}_{h}} \int_{K} \nabla u\cdot\nabla v dx-\sum \limits_{K \in {\xi}_{h}} \int_{\partial K} (\nabla u\cdot {\bf n})vds = \sum \limits_{K \in {\xi}_{h}} \int_{K} fvdx+\sum \limits_{ e \in \Gamma_{N}} \int_{e} g^Nvds
$$
Or using the jump definitions ($v\in V_h$ are element-wise discontinuous)
$$
\sum \limits_{K \in {\xi}_{h}} \int_{K} \nabla u\cdot\nabla v dx-\sum \limits_{ e \in \Gamma_{0}\cup\Gamma^{D}} \int_{e} [v\nabla u]ds = \sum \limits_{K \in {\xi}_{h}} \int_{K} fvdx+\sum \limits_{ e \in \Gamma_{N}} \int_{e} g^Nvds
$$
One can easily verify that $[v\nabla u]=\{ \nabla u\}\cdot [v]+ [\nabla u]\cdot\{ v\}$. Then, using also the fact that $[\nabla u]=0$ ($u$ is assumed to be smooth enough so that $\nabla u$ is continuous), we get
$$
\sum \limits_{K \in {\xi}_{h}} \int_{K} \nabla u\cdot\nabla v dx-\sum \limits_{ e \in \Gamma_{0}\cup\Gamma^{D}} \int_{e} \{ \nabla u\}\cdot [v]ds = \sum \limits_{K \in {\xi}_{h}} \int_{K} fvdx+\sum \limits_{ e \in \Gamma_{N}} \int_{e} g^Nvds
$$
Yet, the left hand side is not coercive, even not symmetric. To handle this, noting that $[u]=0$ along the interior edges ($u$ is assumed to be continuous), we reach at
$$
\sum \limits_{K \in {\xi}_{h}} \int_{K} \nabla u\cdot\nabla v dx-\sum \limits_{ e \in \Gamma_{0}\cup\Gamma^{D}} \int_{e} \{ \nabla u\}\cdot [v]ds -\sum \limits_{ e \in \Gamma_{0}} \int_{e} \{ \nabla v\}\cdot [u]ds
$$
$$
\qquad +\sum \limits_{ e \in \Gamma_{0}} \frac{\sigma}{h_{e}}\int_{e} [u]\cdot [v]ds = \sum \limits_{K \in {\xi}_{h}} \int_{K} fvdx+\sum \limits_{ e \in \Gamma_{N}} \int_{e} g^Nvds
$$
where $h_e$ denote the length of the edge $e$ and $\sigma$ is called the penalty parameter, which is a sufficiently large to have the coercivity. Finally, we add to the both sides the edge integrals on the Dirichlet boundary edges (keeping unknown on the left hand side and imposing Dirichlet boundary condition on the right hand side)
$$
\sum \limits_{K \in {\xi}_{h}} \int_{K} \nabla u\cdot\nabla v dx-\sum \limits_{ e \in \Gamma_{0}\cup\Gamma^{D}} \int_{e} \{ \nabla u\}\cdot [v]ds -\sum \limits_{ e \in \Gamma_{0}\cup\Gamma^{D}} \int_{e} \{ \nabla v\}\cdot [u]ds
$$
$$
\qquad +\sum \limits_{ e \in \Gamma_{0}\cup\Gamma^{D}} \frac{\sigma}{h_{e}}\int_{e} [u]\cdot [v]ds = \sum \limits_{K \in {\xi}_{h}} \int_{K} fvdx+\sum \limits_{ e \in \Gamma_{D}} \int_{e} g^D\left( \frac{\sigma}{h_{e}}v-\nabla v\cdot {\bf n}\right)ds+\sum \limits_{ e \in \Gamma_{N}} \int_{e} g^Nvds
$$
which gives the SIPG formulation.

Now, we give the SIPG discretized system to the problem (\ref{1}) combining with the upwind discretization for the convection part: find $u_h\in V_h$ such that
\begin{equation} \label{ds}
a_{h}(u_{h},v_{h})+b_{h}(u_{h}, v_{h})=l_{h}(v_{h}) \qquad \forall v_h\in V_h,
\end{equation}
\begin{subequations}\label{3}
\begin{align}
a_{h}(u_{h}, v_{h})=& \sum \limits_{K \in {\xi}_{h}} \int_{K} \epsilon \nabla u_{h}\cdot\nabla v_{h} dx + \sum \limits_{K \in {\xi}_{h}} \int_{K} ({\bf b} \cdot \nabla u_{h}+\alpha u_h) v_{h} dx\nonumber \\
&-  \sum \limits_{ e \in \Gamma_{0}\cup\Gamma^{D}} \int_{e} \{\epsilon \nabla v_{h}\} \cdot [u_{h}] ds -\sum \limits_{ e \in \Gamma_{0}\cup\Gamma^{D}}\int_{e} \{\epsilon \nabla u_{h}\} \cdot [v_{h}] ds \nonumber \\
&+ \sum \limits_{K \in {\xi}_{h}}\int_{\partial K^-\setminus\partial\Omega } {\bf b}\cdot {\bf n} (u_{h}^{out}-u_{h}^{in})  v_{h} ds - \sum \limits_{K \in {\xi}_{h}} \int_{\partial K^-\cap \Gamma^{-}} {\bf b}\cdot {\bf n} u_{h}^{in} v_{h}  ds \nonumber  \\
&+ \sum \limits_{ e \in \Gamma_{0}\cup\Gamma^{D}}\frac{\sigma \epsilon}{h_{e}} \int_{e} [u_{h}] \cdot [v_{h}] ds, \nonumber  \\
b_{h}(u_{h}, v_{h}) =& \sum \limits_{K \in {\xi}_{h}} \int_{K} r(u_{h}) v_{h} dx, \nonumber  \\
l_{h}( v_{h})=&  \sum \limits_{K \in {\xi}_{h}} \int_{K} f v_{h} dx
+ \sum \limits_{e \in \Gamma^{D}} \int_e g^D \left( \frac{\sigma \epsilon}{h_{e}} v_{h} -  {\epsilon\nabla v_{h}} \cdot {\bf n} \right) ds  \nonumber \\
&- \sum \limits_{K \in {\xi}_{h}}\int_{\partial K^-\cap \Gamma^{-}} {\bf b}\cdot {\bf n} g^D v_{h}  ds + \sum \limits_{e \in \Gamma^{N}} \int_e g^N v_{h} ds, \nonumber
\end{align}
\end{subequations}
where $\partial K^-$ and $\Gamma^{-}$ indicates the corresponding inflow parts, and $u_{h}^{out}$, $u_{h}^{in}$ denotes the values on an edge from outside and inside of an element $K$, respectively.
The parameter $\sigma\in\mathbb{R}_0^+$ is called the penalty parameter which should be sufficiently large; independent of the mesh size $h$ and the diffusion coefficient $\epsilon$ \citep{riviere08dgm} [Sec. 2.7.1].
We choose the penalty parameter $\sigma$ for the SIPG method depending on the polynomial degree $k$ as
$\sigma=3k(k+1)$ on interior edges and $\sigma=6k(k+1)$ on boundary edges.

\section{Adaptivity}
Most of the convection dominated problems lead to internal/boundary layers and one has to find  accurate approximations in order to handle the nonphysical oscillations. A naive approach is to refine the mesh uniformly. But it is  not desirable as it highly increases the degrees of freedom and refines the mesh unnecessarily in regions where the solutions are smooth.  Instead, the mesh is refined locally using an adaptive strategy. In this section, we describe the adaptive strategy for non-linear diffusion-convection-reaction problems.

\subsection{The adaptive procedure}

\begin{figure}[hbt]
\centering
\setlength{\unitlength}{0.6mm}
\begin{picture}(140,105)

\put(70,100){\oval(20,8)}
\put(62,98){\text{Begin}}
\put(70,96){\vector(0,-1){6}}
\put(40,90){\line(1,0){78}}
\put(40,80){\line(1,0){78}}
\put(40,80){\line(0,1){10}}
\put(118,80){\line(0,1){10}}
\put(42,84){\text{Initialization:mesh, $0<tol, \theta$}}
\put(70,80){\vector(0,-1){10}}
\put(60,70){\line(1,0){24}}
\put(60,63){\line(1,0){24}}
\put(60,63){\line(0,1){7}}
\put(84,63){\line(0,1){7}}
\put(62,65){\text{SOLVE}}
\put(70,63){\vector(0,-1){6}}
\put(45,57){\line(1,0){65}}
\put(45,50){\line(1,0){65}}
\put(45,50){\line(0,1){7}}
\put(110,50){\line(0,1){7}}
\put(46,52){\text{ESTIMATE: compute $\eta$}}
\put(70,50){\vector(0,-1){4}}
\put(72,42){\oval(25,8)}
\put(63,41){\text{$\eta< tol$}}
\put(70,38){\vector(0,-1){6}}
\put(60,33){\text{No}}
\put(45,32){\line(1,0){65}}
\put(45,24){\line(1,0){65}}
\put(45,24){\line(0,1){8}}
\put(110,24){\line(0,1){8}}
\put(70,24){\vector(0,-1){6}}
\put(47,26){\text{MARK: find subset $M_{K}$}}
\put(37,18){\line(1,0){90}}
\put(37,10){\line(1,0){90}}
\put(37,10){\line(0,1){8}}
\put(127,10){\line(0,1){8}}
\put(70,10){\vector(0,-1){6}}
\put(38,12){\text{REFINE: refine triangles $K \in M_{K}$}}
\put(70,0){\oval(20,8)}
\put(65,-2){\text{End}}
\put(32,8){\line(1,0){38}}
\put(32,8){\line(0,1){68}}
\put(32,76){\vector(1,0){38}}
\put(135,20){\text{Yes}}
\put(133,0){\vector(-1,0){53}}
\put(133,0){\line(0,1){43}}
\put(84,43){\line(1,0){49}}
\end{picture}
\caption{Adaptive strategy}
\label{adapstr}
\end{figure}
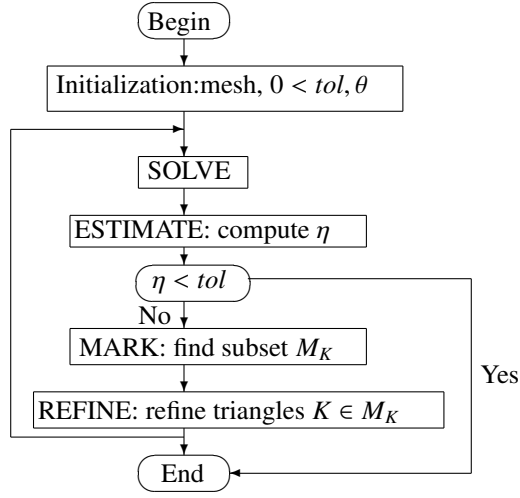

Our adaptive algorithm is based on the standard adaptive finite element (AFEM) iterative loop (Fig.\ref{adapstr}).
The first step, SOLVE, is to solve the SIPG discretized system (\ref{ds}) on a given triangulation $\xi_h$. The ESTIMATE step is the key part of the adaptive procedure, by which we are able to determine the elements having large error to be refined using computed solution and given data (a posteriori). As an estimator, we use a residual based error indicator based on the modification of the error estimator given in Sch\"{o}tzau and Zhu \citep{schotzau09rae} for a single linear convection dominated diffusion-convection-reaction equation  to the diffusion-convection equation with non-linear reaction mechanism, which is robust, i.e. independent of the P\'{e}clet number. To do this, we include in the a posteriori error indicator the non-linear reaction term as local contributions to the cell residuals and not to the interior/boundary edge residuals [Chp. 5.1.4, \citep{verfurth13aee}]. Let $u_h$ be the solution to (\ref{ds}). Then, for each element $K \in {\xi}_{h}$, we define the local error indicators $\eta_K^2$ as

\begin{eqnarray} \label{res}
\eta_K^2= \eta_{R_K}^2  + \eta_{E_K^0}^2 + \eta_{E_K^D}^2+ \eta_{E_K^N}^2,
\end{eqnarray}
In (\ref{res}), $\eta_{R_K}$ denote the cell residual
\begin{eqnarray}
\eta_{R_K}^2&=&\rho_K^2\| f-\alpha u_h+\epsilon\Delta u_h-{\bf b}\cdot\nabla u_h-r(u_h)\|_{L^2(K)}^2, \nonumber
\end{eqnarray}
while, $\eta_{E_K^0}$, $\eta_{E_K^D}$ and $\eta_{E_K^N}$ stand for the edge residuals coming from the jump of the numerical solution on the interior, Dirichlet boundary and Neumann boundary edges, respectively
\begin{eqnarray}
\eta_{E_K^0}^2 &=& \sum \limits_{e \in \partial K\cap\Gamma_{0}}\left(\frac{1}{2}\epsilon^{-\frac{1}{2}}\rho_e\| [\epsilon\nabla u_h]\|_{L^2(e)}^2+\frac{1}{2}(\frac{\epsilon\sigma}{h_e}+\kappa h_e+\frac{h_e}{\epsilon})\| [u_h]\|_{L^2(e)}^2\right), \nonumber \\
\eta_{E_K^D}^2 &=& \sum \limits_{e \in \partial K\cap\Gamma^{D}}(\frac{\epsilon\sigma}{h_e}+\kappa h_e+\frac{h_e}{\epsilon})\| g^D-u_h\|_{L^2(e)}^2, \nonumber \\
\eta_{E_K^N}^2 &=& \sum \limits_{e \in \partial K\cap\Gamma^{N}}\epsilon^{-\frac{1}{2}}\rho_e\| g^N-\epsilon\nabla u_h\cdot {\bf n}\|_{L^2(e)}^2, \nonumber
\end{eqnarray}
as in Sch\"{o}tzau and Zhu \citep{schotzau09rae} with the modified \citep{hoppe08caa,schotzau09rae} extra term $\eta_{E_K^N}^2$ corresponding to the local indicator on the Neumann boundary edges. The weights $\rho_K$ and $\rho_e$, on an element $K$, are defined as
$$
\rho_{K}=\min\{h_{K}\epsilon^{-\frac{1}{2}}, \kappa^{-\frac{1}{2}}\}, \; \rho_{e}=\min\{h_{e}\epsilon^{-\frac{1}{2}}, \kappa^{-\frac{1}{2}}\},
$$
for $\kappa \neq 0$. When $\kappa =0$, we take $\rho_{K}=h_{K}\epsilon^{-\frac{1}{2}}$ and $\rho_{e}=h_{e}\epsilon^{-\frac{1}{2}}$. Then, our a posteriori error indicator is given by
$$
\eta=\left( \sum \limits_{K\in{\xi}_{h}}\eta_K^2\right)^{1/2}.
$$
We also introduce the data approximation error,
$$
\Theta^2=\Theta^2(f)+\Theta^2(u^D)+\Theta^2(u^N)
$$
where

\begin{eqnarray}
\Theta^2(f) &=&\sum \limits_{K \in \xi}\rho_K^2(\| f-f_h\|_{L^2(K)}^2 + \| ({\bf b} -{\bf b}_h)\cdot\nabla u_h\|_{L^2(K)}^2 +\| (\alpha -\alpha_h)u_h\|_{L^2(K)}^2), \nonumber \\
\Theta^2(u^D) &=& \sum \limits_{e \in \Gamma^D} (\frac{\epsilon\sigma}{h_e}+\kappa h_e+\frac{h_e}{\epsilon})\| g^D-\hat{g}^D\|_{L^2(e)}^2, \nonumber \\
\Theta^2(u^N) &=& \sum \limits_{e \in \Gamma^N} \epsilon^{-\frac{1}{2}}\rho_e\| g^N-\hat{g}^N\|_{L^2(e)}^2, \nonumber
\end{eqnarray}
according to \citep{schotzau09rae}, with $\hat{g}^D$ and $\hat{g}^N$ denoting the mean integrals of $g^D$ and $g^N$, respectively.

In the MARK step, if the given tolerance is not satisfied, we determine the set of elements $M_K\subset {\xi}_{h}$ to be refined using the error indicator defined in (\ref{res}). To do this, we use the bulk criterion proposed by D\"{ofler} \citep{dorfler96caa}, by which the approximation error is decreased by a fixed factor for each loop. In the light of bulk criterion, we choose the set of elements $M_K\subset {\xi}_{h}$ satisfying
$$
\sum \limits_{K \in M_K}\eta_K^2 \geq  \theta \sum \limits_{K \in \xi_h}\eta_K^2,
$$
for a user defined parameter $0<\theta<1$. Here, bigger $\theta$ results in more refinement of triangles in a single loop, where, smaller $\theta$ causes more refinement loops.

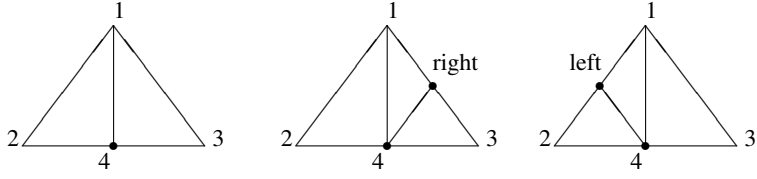
\begin{figure}[hbt]
\centering

\setlength{\unitlength}{2mm}
\begin{picture}(50, 11)

\put(2,2){\line(1,0){12}}
\put(2,2){\line(3,4){6}}
\put(8,2){\line(0,1){8}}
\put(14,2){\line(-3,4){6}}
\put(1,2){\small 2 }
\put(14.5,2){\small 3 }
\put(8,10.5){\small 1 }
\put(7,0.5){\small 4 }
\put(8,2){\circle*{0.6} }

\put(20,2){\line(1,0){12}}
\put(20,2){\line(3,4){6}}
\put(26,2){\line(0,1){8}}
\put(32,2){\line(-3,4){6}}
\put(19,2){\small 2 }
\put(32.5,2){\small 3 }
\put(26,10.5){\small 1 }
\put(25,0.5){\small 4 }
\put(26,2){\circle*{0.55} }
\put(26,2){\line(3,4){3}}
\put(29,6){\circle*{0.55} }
\put(29,7){\small right }

\put(37,2){\line(1,0){12}}
\put(37,2){\line(3,4){6}}
\put(43,2){\line(0,1){8}}
\put(49,2){\line(-3,4){6}}
\put(36,2){\small 2 }
\put(49.5,2){\small 3 }
\put(43,10.5){\small 1 }
\put(42,0.5){\small 4 }
\put(43,2){\circle*{0.55} }
\put(43,2){\line(-3,4){3}}
\put(40,6){\circle*{0.55} }
\put(38,7){\small left }

\end{picture}

\caption{Bisection of a triangle}.
\label{bisect}
\end{figure}

Finally, REFINE step, we refine the marked elements $K\in M_K$ using the newest vertex bisection method \citep{chen08fem}. This process can be summarized as (see Fig.\ref{bisect}): for each element $K\in\xi_h$, we label one vertex of $K$ as a newest vertex. The opposite edge of the newest vertex is called as the refinement edge. Then, a triangle is bisected to two new children triangles by connecting the newest vertex to the midpoint of the refinement edge, and this new vertex created at the midpoint of the refinement edge is assigned to be the newest vertex of the children. Following a similar rule, these two children triangles are bisected to obtain four children elements belonging to the father element (the refined triangle $K\in M_K$). After bisecting all $K\in M_K$, we also divide some elements $K\in\xi_h\setminus M_K$ to keep the conformity of the mesh, i.e. hanging nodes are not allowed.

In the case of coupled problems, instead of a single component problem, we refine the elements being the union of the set of the elements to be refined for each component, i.e.,  let $\eta_K^1$ and $\eta_K^2$ be the computed local error indicators corresponding to each unknown component of a two component system. Next, we determine the set of elements $M_K^1$ and $M_K^2$ satisfying
$$
\sum \limits_{K \in M_K^1}(\eta_K^1)^2 \geq \theta \sum \limits_{K \in \xi_h}(\eta_K^1)^2   \; , \qquad \sum \limits_{K \in M_K^2}(\eta_K^2)^2 \geq \theta \sum \limits_{K \in \xi_h}(\eta_K^2)^2.
$$
Then, we refine the marked elements $K\in M_K^1\cup M_K^2$ using the newest vertex bisection method. The adaptive procedure ends after a sequence of mesh refinements up to attain a solution with an estimated error within a prescribed tolerance. Numerical studies show the capability of the error indicator to find the layers properly.

\subsection{Reliability and efficiency of a posteriori error estimator}

In order to measure the error, we use the energy norm

$$
||| v |||^2=\sum \limits_{K \in {\xi}_{h}}(\| \epsilon\nabla v\|_{L^2(K)}^2 +\kappa \| v\|_{L^2(K)}^2) + \sum \limits_{e \in \Gamma_{0}\cup\Gamma^{D}}\frac{\epsilon\sigma}{h_{e}}\| [v]\|_{L^2(e)}^2,
$$
and the semi-norm
\begin{equation} \label{semin}
|v|_C^2=|{\bf b}v|_*^2+\sum\limits_{e \in \Gamma_{0}\cup\Gamma^{D}}(\kappa h_e+ \frac{h_e}{\epsilon})\| [v]\|_{L^2(e)}^2,
\end{equation}
where
$$
|u|_*=\mathop{\text{sup}}_{w\in H_0^1(\Omega )\setminus \{ 0\}}\frac{\int_{\Omega}u\cdot \nabla w dx}{||| w|||},
$$
and the constant $\kappa \geq 0$ satisfies
\begin{equation} \label{erassmp}
\alpha-\frac{1}{2}\nabla\cdot{\bf b} (x) \geq \kappa, \qquad \| -\nabla\cdot{\bf b} +\alpha\|_{L^{\infty}(\Omega)}\leq \kappa^*\kappa,
\end{equation}
for a non-negative $\kappa^*$. The first condition in (\ref{erassmp}) is needed to make the bilinear form $a_h$ coercive (well-posedness of the linear part of the problem), and the latter is used to prove the reliability of our a posteriori error estimator \citep{schotzau09rae}. The terms $|{\bf b} v|_*^2$ and $\frac{h_e}{\epsilon}\| [v]\|_{L^2(e)}^2$ in (\ref{semin}) are used to bound the convection part, whereas, the term $\kappa h_e\| [v]\|_{L^2(e)}^2$ is used to bound the linear reaction part of the discrete system. In order to bound the non-linear reaction part, we use the boundedness property (2a) [Chp. 5.1.1-4, \citep{verfurth13aee}].
Then, for the solution $u$ to the scalar equation of (\ref{1}), following the procedures in \citep{schotzau09rae} and using the boundedness of the non-linear reaction term, we can easily obtain the a posteriori error bounds

\begin{eqnarray}
||| u-u_h||| + | u-u_h|_C \lesssim \eta + \Theta  \qquad &(\text{reliability}), \nonumber \\[0.2cm]
\eta \lesssim ||| u-u_h||| + | u-u_h|_C + \Theta \qquad &(\text{efficiency}). \nonumber
\end{eqnarray}

\section{Comparison with the Galerkin least squares finite element method (GLSFEM)}
For linear PDEs, the weak form in the standard Galerkin method is obtained by multiplying the differential equation with a test function $v$ and integrating over a suitable
function space $V$
$$
(\mathcal{L}u,v) = (f,v), \quad \forall v \in V
$$
where $\mathcal{L}=-\epsilon\Delta +{\bf b}\cdot\nabla+\alpha$ is the linear part of the diffusion-convection-reaction equation (1). Defining the residual as $R(u)=f-\mathcal{L}u$, the standard Galerkin method can be interpreted in form of the residual orthogonality $(R(u),v)=0$. In the case of non self-adjoint
differential operators like the the diffusion-convection-reaction operator $\mathcal{L}$, it can happen that $(\mathcal{L}u,v)$ is not coercive or symmetric on $V$, and the resulting FEM discretization may be unstable.

For transport problems,  another popular approach is based on the least squares formulation of the Galerkin FEM.
Let us write the model problem (1) as
\begin{subequations}\label{model2}
\begin{align}
\mathcal{L}u +r(u) &= f \quad  \text{ in } \; \Omega\\
u &= g  \quad \; \text{on } \; \Gamma,
\end{align}
\end{subequations}
Define the least-squares functional
$$
J(u):=\frac{1}{2}\| \mathcal{L}u+r(u)-f\|^2_{L^2(\Omega)}.
$$
A minimizer of $J(u)$ is obtained by
$$
\lim_{t\mapsto 0}\frac{d}{dt}J(u+tv)=0 \; , \quad \forall v \in V
$$
which yields the least-squares term
$$
\tilde{J}_{\Omega}(u,v):= (\mathcal{L}u+r(u)-f,\mathcal{L}v+r'(u)v)_{L^2(\Omega)}.
$$
For linear problems with $r(u)=0$, the least squares Galerkin method reduces to the minimization problem
$$
F(u) = \min_{v\in V} F(v)
$$
where the functional $F(\cdot)$ is defined by
$$
F(v) = \frac{1}{2} \mid\mid  \mathcal{L}v - f\mid \mid ^2_{L^2(\Omega)} .
$$
The first order optimality condition leads to the least squares Galerkin method
$$
(\mathcal{L}u,\mathcal{L}v)  = (f,\mathcal{L}v)\; , \quad \forall v \in V.
$$
The bilinear form $(\mathcal{L}u,\mathcal{L}v)$ is symmetric and coercive and has stronger stability properties compared to the standard Galerkin method.

There are many publications on the Galerkin least squares finite element methods (GLSFEM). We mention here two books \citep{bochev09lsfem,jiang98lsfem} and the review article
\citep{bochev98femlserev}. There are mainly two variants of the GLSFEMs; the stabilized and the direct versions.\\

\noindent {\it Stabilized finite elements method \citep{hughes89fem}}: The standard (continuous) Galerkin FEM for the problem (\ref{model2}) reads: find $u_h\in U_h\subset U$ such that
\begin{eqnarray}\label{stG}
a(u_h,v_h) +(r(u_h),v_h)_{L^2(\Omega)} &=& (f,v_h)_{L^2(\Omega)} \; , \quad \, \forall v_h\in V_h\subset V
\end{eqnarray}
where $a(u,v)=(\epsilon\nabla u+{\bf b}\cdot\nabla u+\alpha u,v)_{L^2(\Omega)}$ is the standard bilinear form to the linear part of (\ref{model2}). The general stabilized FEMs formulation reads as: for all $v_h\in V_h\subset V$, find $u_h\in U_h\subset U$ such that
\begin{eqnarray}\label{stFEM}
a(u_h,v_h) +(r(u_h),v_h)_{L^2(\Omega)} + {\bf \sum_K \tau_K  S_K(u_h,v_h)} = (f,v_h)_{L^2(\Omega)}
\end{eqnarray}
where the stabilization parameter is defined on each element $K$ as \citep{hauke02subgrid}
$$
\tau_K = \frac{1}{\frac{4\epsilon }{h^2}+\frac{2|{\bf b}|}{h}+|\alpha |}.
$$
One way to proceed GLSFEMs is then to use the least-squares term $\tilde{J}_{K}(u,v)$ as the stabilization term $S_K$ in (\ref{stFEM}), i.e.: for all $v_h\in V_h$, find $u_h\in U_h$ such that
\begin{eqnarray}\label{stLSFEM}
a(u_h,v_h) +(r(u_h),v_h)_{L^2(\Omega)} + \sum_K \tau_K  \tilde{J}_{K}(u_h,v_h) = (f,v_h)_{L^2(\Omega)}
\end{eqnarray}
Note that, being another stabilized FEM, streamline upwind Petrov-Galerkin (SUPG) method is obtained by setting
$$
S_K(u_h,v_h)=(\mathcal{L}u_h+r(u_h)-f, {\bf b}\cdot\nabla v_h)_{L^2(K)}
$$
with different choices of the parameter $\tau_K$.\\

\noindent {\it The direct variant of  GLSFEM:} The second way to proceed the GLSFEMs is to consider and discretize just the least-squares term $\tilde{J}_{\Omega}(u,v)$. One may solve this problem in a straightforward manner: for all $v_h\in V_h\subset H^2(\Omega)\cap V$, find $u_h\in U_h\subset H^2(\Omega)\cap U$ such that $\tilde{J}_{\Omega}(u_h,v_h) = 0$, i.e.
\begin{eqnarray*}
\int_{\Omega}(\mathcal{L}u_h+r(u_h))(\mathcal{L}v_h+r'(u_h)v_h)dx &=& \int_{\Omega}f(\mathcal{L}v_h+r'(u_h)v_h)dx
\end{eqnarray*}
which is not only a fourth order problem but also the solution and trial subspaces $U_h$ and $V_h$ need to consist of continuously differentiable functions making it complicated to construct bases functions (standard finite element spaces cannot be used anymore) and the assembly of the stiffness matrix. The condition number of the stiffness matrix is order of
$\mathcal{O}(h^{-4})$ instead of order $\mathcal{O}(h^{-2})$ for the standard Galerkin FEM. Hence, this approach is impractical. Instead, being the most common practical way, the problem (\ref{model2}) is converted into a first-order system as \citep{houston02hplsfem,bochev98femlserev}.
\begin{eqnarray*}
p - \nabla u &=& 0 \quad \, \text{ in } \; \Omega\\
-\epsilon\nabla \cdot p +{\bf b}\cdot\nabla u +\alpha u + r(u) &=& f \quad \, \text{ in } \; \Omega\\
u &=& g  \quad \text{ on } \; \Gamma
\end{eqnarray*}
Then, we define now the least-square functional for $z=(p,u)^T$ as
$$
J(z):=\frac{1}{2}\| p - \nabla u\|^2_{L^2(\Omega)} + \frac{1}{2}\| -\epsilon\nabla \cdot p +{\bf b}\cdot\nabla u +\alpha u + r(u) - f\|^2_{L^2(\Omega)}
$$
A minimizer of $J(z)$ is obtained by the identity
$$
\lim_{t\mapsto 0}\frac{d}{dt}J(z+tv)=0 \; , \quad \forall v
$$
which yields a least-squares term of order two. Using this approach, we solve the resulting least-squares term, which is a second-order problem now, using (discontinuous) finite elements solution and trial spaces $S_h\subset H^1(\Omega,div)\times U$ ($S_h\subset H^1(\Omega,div)\times H^1(\Omega)$) and $T_h\subset H^1(\Omega ,div)\times V$ ($T_h\subset H^1(\Omega,div)\times H^1(\Omega )$), respectively. The condition number of the stiffness matrix is retained as $\mathcal{O} (h^{-2})$ as in the standard Galerkin method \citep{bochev98femlserev}. For convection dominated problems, the resulting linear systems of equations are solved usually with preconditioned conjugate gradient method due to large condition numbers, as reported in \citep{lazarov00lsd} for GLSFEM solution of singularly perturbed diffusion-convection problems.

In order to compare the GLSFEM with the DGFEM, we consider the linear problem \citep{yucel13doc}
\begin{equation} \label{lconvdiff}
-\epsilon\Delta u +{\bf b}\cdot\nabla u+\alpha u = f \quad \, \hbox{ in } \; (0,1)^2
\end{equation}
with $\epsilon = 10^{-6}$, $\beta =(2,3)^T$ and $\alpha =1$. The load function $f$ and Dirichlet boundary conditions are chosen so that the exact solution is
$$
u(x_1,x_2)=\frac{\pi }{2}\arctan \left( \frac{1}{\sqrt{\epsilon}}(-0.5x_1+x_2-0.25) \right)
$$

\begin{figure}[htb]
\centering
\subfloat{\includegraphics[width=0.8\textwidth]{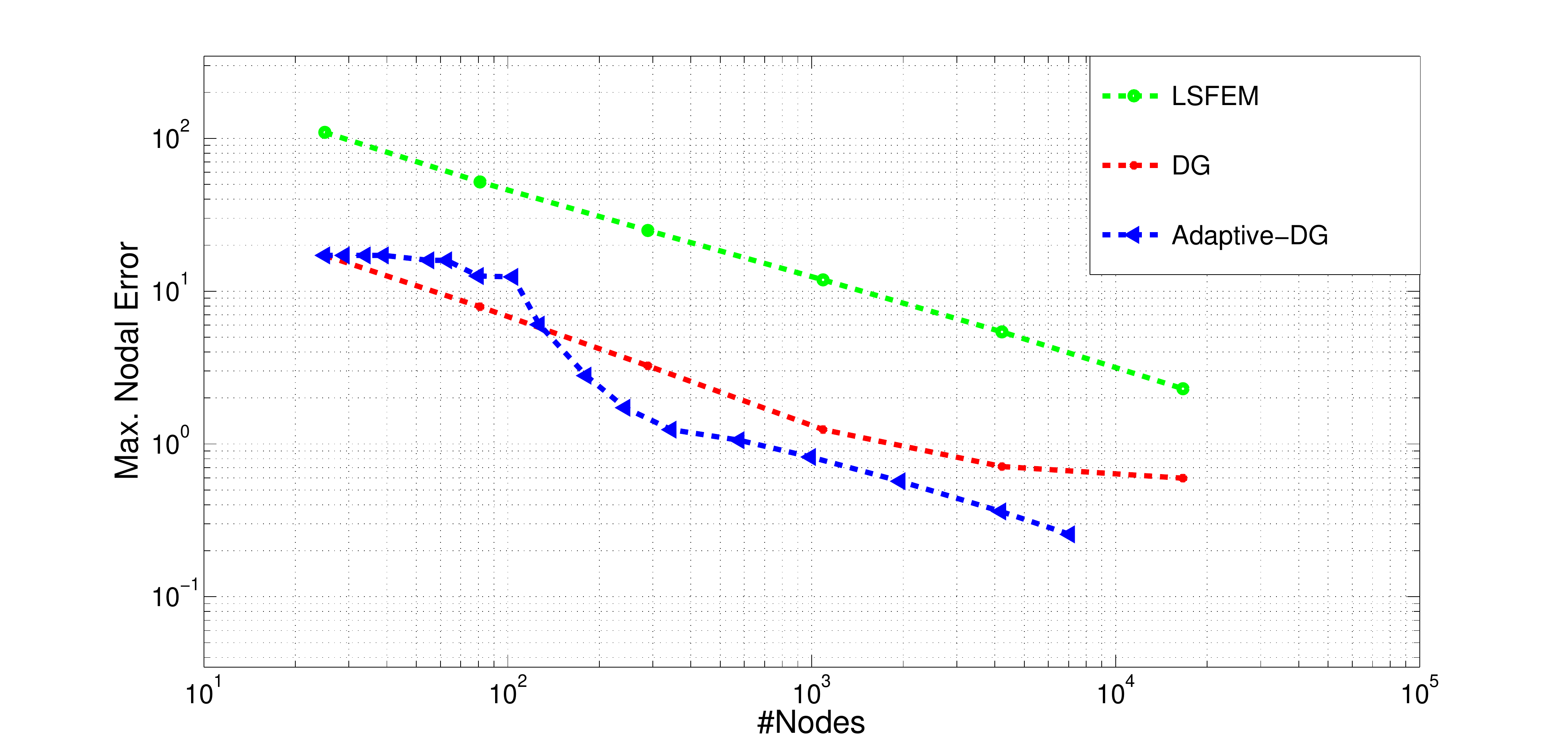}}
\caption{Comparion of the GLSFEM and DGFEM for linear convection dominated problem (\ref{lconvdiff})  }
\label{glserr}
\end{figure}

As we  can see from the from the Fig.\ref{glserr}, the DGFEMs produce smaller errors than the GLSFEM for the convection dominated problem (\ref{lconvdiff}).

\section{Efficient solution of linear systems}\label{linear}

The approximate solution to the discrete problem (\ref{ds}) has the form
$$
u_h=\sum_{i=1}^{Nel}\sum_{l=1}^{Nloc}U_l^i\phi_l^i
$$
where $\phi_l^i$'s are the basis polynomials spanning the DGFEM space $V_h$, $U_l^i$'s are the unknown coefficients to be found, $Nel$ denotes the number of triangles and $Nloc$ is the number of local dimension depending on the degree of polynomials $k$, for instance, for $k=2$ we have $Nloc=6$ (in 2D, $Nloc=(k+1)(k+2)/2$). In DG methods, we choose the piecewise basis polynomials $\phi_l^i$'s in such a way that each basis function has only one triangle as a support, i.e. we choose on a specific triangle $K_e$, $e\in\{ 1,2,\ldots , Nel\}$, the basis polynomials $\phi_l^e$ which are zero outside the triangle $K_e$, $l=1,2,\ldots , Nloc$. By this construction, the stiffness matrix obtained by DG methods has a block structure, each of which related to a triangle (there is no overlapping as in continuous FEM case). The product $dof:=Nel*Nloc$ gives the degree of freedom in DG methods. Inserting the linear combination of $u_h$ in (\ref{ds}) and choosing the test functions as $v_h=\phi_l^i$, $l=1,2,\ldots , Nloc$, $i=1,2,\ldots , Nel$, the discrete residual of the system (\ref{ds}) in matrix vector form is given by
$$
R({\bf U})=S{\bf U}+h({\bf U})-L
$$
where ${\bf U}\in\mathbb{R}^{dof}$ is the vector of unknown coefficients $U_l^i$'s, $S\in\mathbb{R}^{dof\times dof}$ is the stiffness matrix corresponding to the bilinear form $a_h(u_h,v_h)$, $h\in\mathbb{R}^{dof}$ is the vector function of $U$ related to the non-linear form $b_h(u_h,v_h)$ and $L\in\mathbb{R}^{dof}$ is the vector to the linear form $l_h(v_h)$. The explicit definitions are given by
$$
S=
\begin{bmatrix}
S_{11} & S_{12} & \cdots & S_{1,Nel} \\
S_{21} & S_{22} &  & \vdots \\
\vdots &  & \ddots &  \\
S_{Nel,1} & \cdots &  & S_{Nel,Nel}
\end{bmatrix} \; , \quad {\bf U}=
\begin{bmatrix}
 {\bf U}_1 \\
{\bf U}_2 \\
\vdots \\
{\bf U}_{Nel}
\end{bmatrix}
$$

$$
h({\bf U})=
\begin{bmatrix}
 {\bf h}_1 \\
{\bf h}_2 \\
\vdots \\
{\bf h}_{Nel}
\end{bmatrix} \; , \quad L=
\begin{bmatrix}
 {\bf L}_1 \\
{\bf L}_2 \\
\vdots \\
{\bf L}_{Nel}
\end{bmatrix}
$$
where all the block matrices have dimension $Nloc$:
$$
S_{ji}=
\begin{bmatrix}
a_h(\phi_1^i,\phi_1^j) & a_h(\phi_2^i,\phi_1^j) & \cdots & a_h(\phi_{Nloc}^i,\phi_1^j) \\
a_h(\phi_1^i,\phi_2^j)  & a_h(\phi_2^i,\phi_2^j)  &  & \vdots \\
\vdots &  & \ddots &  \\
a_h(\phi_1^i,\phi_{Nloc}^j)  & \cdots &  & a_h(\phi_{Nloc}^i,\phi_{Nloc}^j)
\end{bmatrix} \; , \quad {\bf U}_i=
\begin{bmatrix}
 U_1^i \\
U_2^i \\
\vdots \\
U_{Nloc}^i
\end{bmatrix}
$$

$$
{\bf h}_i=
\begin{bmatrix}
 b_h(u_h,\phi_1^i) \\
b_h(u_h,\phi_2^i) \\
\vdots \\
b_h(u_h,\phi_{Nloc}^i)
\end{bmatrix} \; , \quad {\bf L}_i=
\begin{bmatrix}
 l_h(\phi_1^i) \\
l_h(\phi_2^i) \\
\vdots \\
l_h(\phi_{Nloc}^i)
\end{bmatrix}
$$

\begin{figure}[hbt]
\centering
\setlength{\unitlength}{2mm}
\begin{picture}(55, 17)
\put(15,2){\line(-3,4){8}}
\put(7,2){\vector(0,1){13}}
\put(7,2){\vector(1,0){10}}
\put(9,5){\text{$\hat{T}$}}
\put(7.5,15){\small$x_2$ }
\put(16,3){\small$x_1$ }
\put(4,0){\small$(0,0)$ }
\put(13,0){\small$(1,0)$ }
\put(3,12){\small$(0,1)$ }

\put(22,15){\small$x_1=\frac{(1+z_1)(1-z_2)}{4}$ }
\put(22,12){\small$x_2=\frac{1+z_2}{2}$ }
\put(34,9){\vector(-1,0){12}}
\put(22,6){\vector(1,0){12}}
\put(22,3){\small$z_1=\frac{2x_1}{1-x_2}-1$ }
\put(22,0){\small$z_2=2x_2-1$ }

\put(42,2){\line(1,0){10}}
\put(42,2){\line(0,1){10}}
\put(42,12){\line(1,0){10}}
\put(52,2){\line(0,1){10}}
\put(47,7){\vector(0,1){7}}
\put(47,7){\vector(1,0){7}}
\put(44,5){\text{$\hat{Q}$}}
\put(47.5,14){\small$z_2$ }
\put(54,8){\small$z_1$ }
\put(39,0){\small$(-1,-1)$ }
\put(50,0){\small$(1,-1)$ }
\put(39,13){\small$(-1,1)$ }
\put(50.5,13){\small$(1,1)$ }
\end{picture}
\caption{Collapsed coordinate transform between reference triangle and reference square}.
\label{trans}
\end{figure}
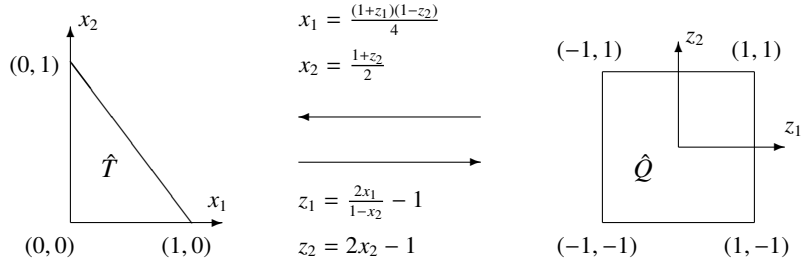

\begin{figure}[htb]
\centering
\includegraphics[width=0.8\textwidth]{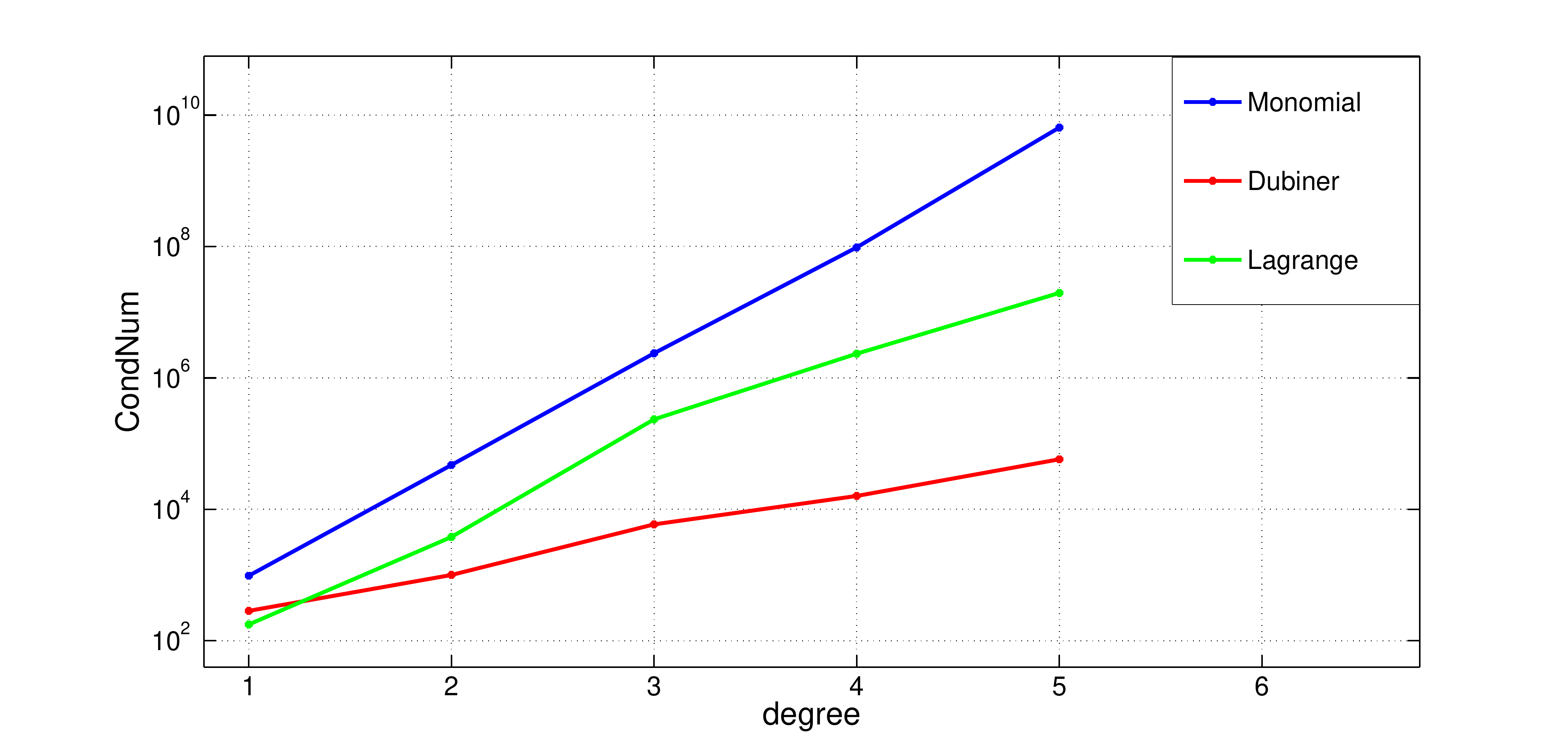}
\caption{Degree vs. condition number of the stiffness matrix: comparison for different type of basis functions for the linear model problem (\ref{1}) with $r(u)=0$   }
\label{basiscond}
\end{figure}

Obviously, the condition number of the stiffness matrix obtained by the SIPG discretization increases by the degree $k$ of basis polynomials. One of the reasonable ways to handle this drawback is to choose a suitable set of basis polynomials. There are a variety of basis polynomial functions such as Lagrange shape functions, monomial bases, Legendre polynomials. In our study, we use the orthogonal Dubiner basis defined on the reference triangle \citep{deng05aao}
$$
\hat{T}=\{ {\bf x}=(x_1,x_2) | \; 0\leq x_1,x_2 \leq 1 \}
$$
(all the integral terms above are computed on this reference triangle using an affine map between the reference triangle and physical triangles). The construction of such basis polynomials based on the collapsed coordinate transform between the reference triangle $\hat{T}$ and the reference square $\hat{Q}=[-1,1]^2$ (see Fig.\ref{trans}).
First, the basis polynomials on the square $\hat{Q}$ is formed by a generalized tensor product of the Jacobi polynomials on the interval $[-1,1]$, and then, these basis polynomials are transformed to the reference triangle $\hat{T}$ using the collapsed coordinate transform in Fig.\ref{trans}. The explicit forms of Dubiner basis polynomials on the reference triangle $\hat{T}$ are given by
\begin{eqnarray*}
\phi_{mn}(x_1,x_2) &=& (1-z_2)^mP_m^{0,0}(z_1)P_n^{2m+1,0}(z_2) \\
&=& 2^m(1-x_2)^mP_m^{0,0}(\frac{2x_1}{1-x_2}-1)P_n^{2m+1,0}(2x_2-1) \; , \quad 0\leq m,n,m+n \leq Nloc
\end{eqnarray*}
where $P_n^{\alpha , \beta }(x)$'s denote the corresponding $n$-th order Jacobi polynomials on the interval $[-1,1]$, which are orthogonal polynomials under the Jacobi weight $(1-x)^{\alpha }(1+x)^{\beta}$, i.e.
$$
\int_{-1}^1 (1-x)^{\alpha }(1+x)^{\beta}P_m^{\alpha , \beta }(x)P_n^{\alpha , \beta }(x)dx=\delta_{mn}
$$
This property of the Jacobi polynomials yields the orthogonality of the Dubiner basis on the reference triangle $\hat{T}$ as
$$
\iint_{\hat{T}} \phi_{mn}(x_1,x_2)\phi_{ij}(x_1,x_2)dx_1dx_2=\frac{1}{8}\delta_{mi}\delta_{nj}.
$$
The advantage of the Dubiner basis is that its orthogonality leads to diagonal mass matrix by which one may obtain better-conditioned matrices compared to the other basis polynomials (see Fig.\ref{basiscond}), and it provides high accuracy in the approximation of the integrals.

\subsection{Effect of the penalty parameter}

\begin{figure}[htb]
\centering
\includegraphics[width=0.8\textwidth]{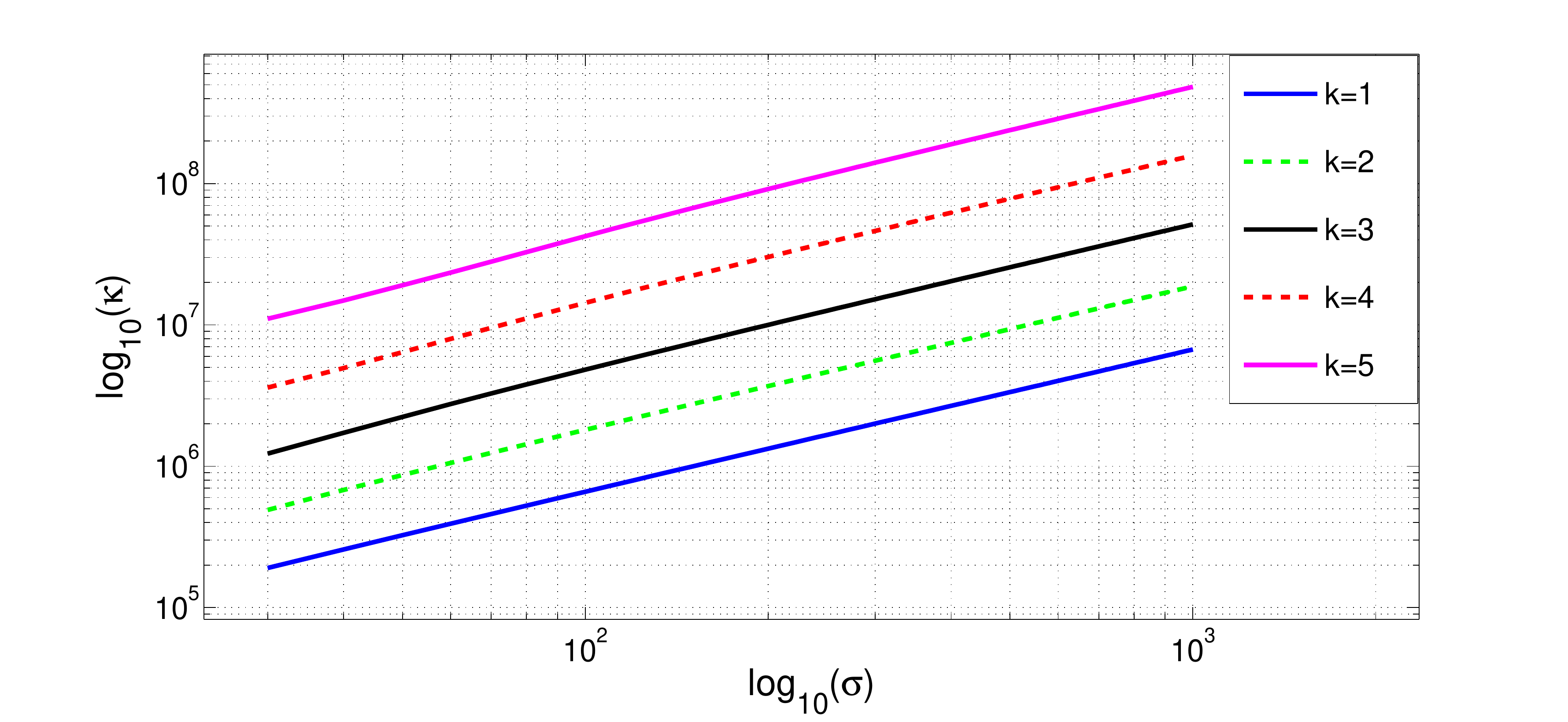}
\caption{Condition number of the stiffness matrix of the SIPG method as a function of penalty parameter $\sigma$ with different polynomial degree $k$ for the Poisson equation (\ref{poisson2})}
\label{penalty_cond}
\end{figure}

\begin{figure}[htb]
\centering
\includegraphics[width=0.8\textwidth]{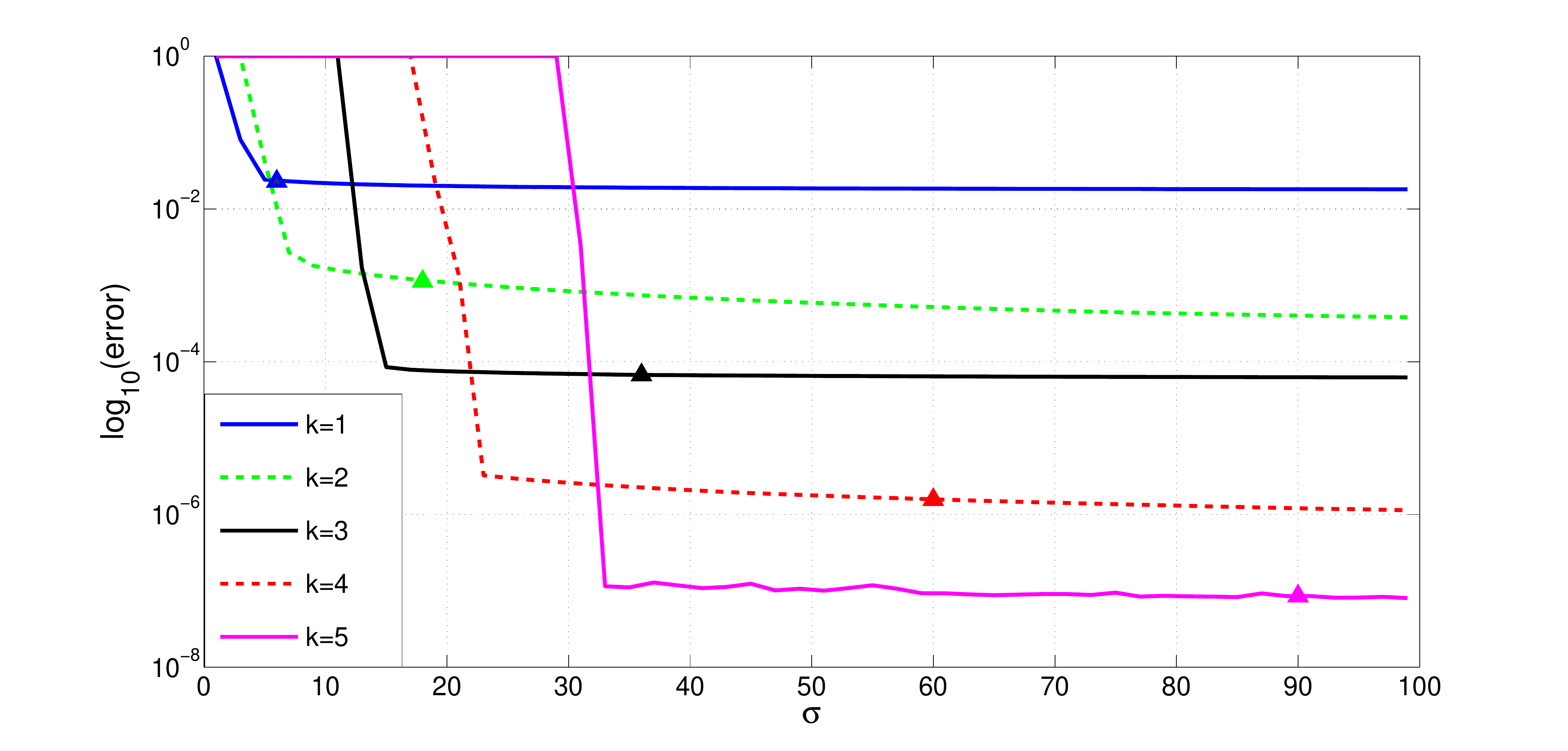}
\caption{Maximum nodal errors of the SIPG approximation as a function of penalty parameter $\sigma$ with different polynomial degree $k$ for the Poisson equation  (\ref{poisson2})}
\label{penalty_err}
\end{figure}

The penalty parameter $\sigma$ in the SIPG formulation (\ref{ds}) should be selected sufficiently large to ensure the coercivity of the bilinear form
[\citep{riviere08dgm}, Sec. 27.1], which is needed for the stability of the convergence of the DG method. It ensures that the matrix arising from the DG discretization of the diffusion part of (\ref{ds}) is symmetric positive definite. At the same time it should not be too large since the conditioning of the matrix arising from the bilinear form increases linearly by the penalty parameter (see Fig.\ref{penalty_cond}). In the literature, several choices of the penalty parameter are suggested. In \citep{epshteyn07epp}, computable lower bounds are derived, and in \citep{dobrev08psi}, the penalty parameter is chosen depending on the diffusion coefficient $\epsilon$. The effect of the penalty parameter on the condition number was discussed in detail for the DG discretization of the Poisson equation in
\citep{castillo12pdg} and in \citep{vuik14fls} for layered reservoirs with strong permeability
contrasts, e.g. $\epsilon$ varying between $10^{-1}$ and $10^{-7}$. Since the penalty parameter, in SIPG formulation, is mainly related to the Laplace operator, to examine the effect of the penalty parameter, we study on the Poisson problem (pure elliptic case )
\begin{eqnarray}\label{poisson2}
-\Delta u &=& f \quad \, \text{ in } \; (0,1)^2
\end{eqnarray}
with the appropriate load function $f$ and Dirichlet boundary conditions using the exact solution $u(x)=\sin (\pi x_1)\sin (\pi x_2)$. In  Fig.\ref{penalty_err}, we have plotted the maximum nodal errors for the Poisson problem (\ref{poisson2}) depending on the penalty parameter to show the instability bound of the scheme for different degrees of bases, where  the triangular symbols indicate our choice $\sigma =3k(k+1)$.

\begin{figure}[htb]
\centering
\includegraphics[width=0.8\textwidth]{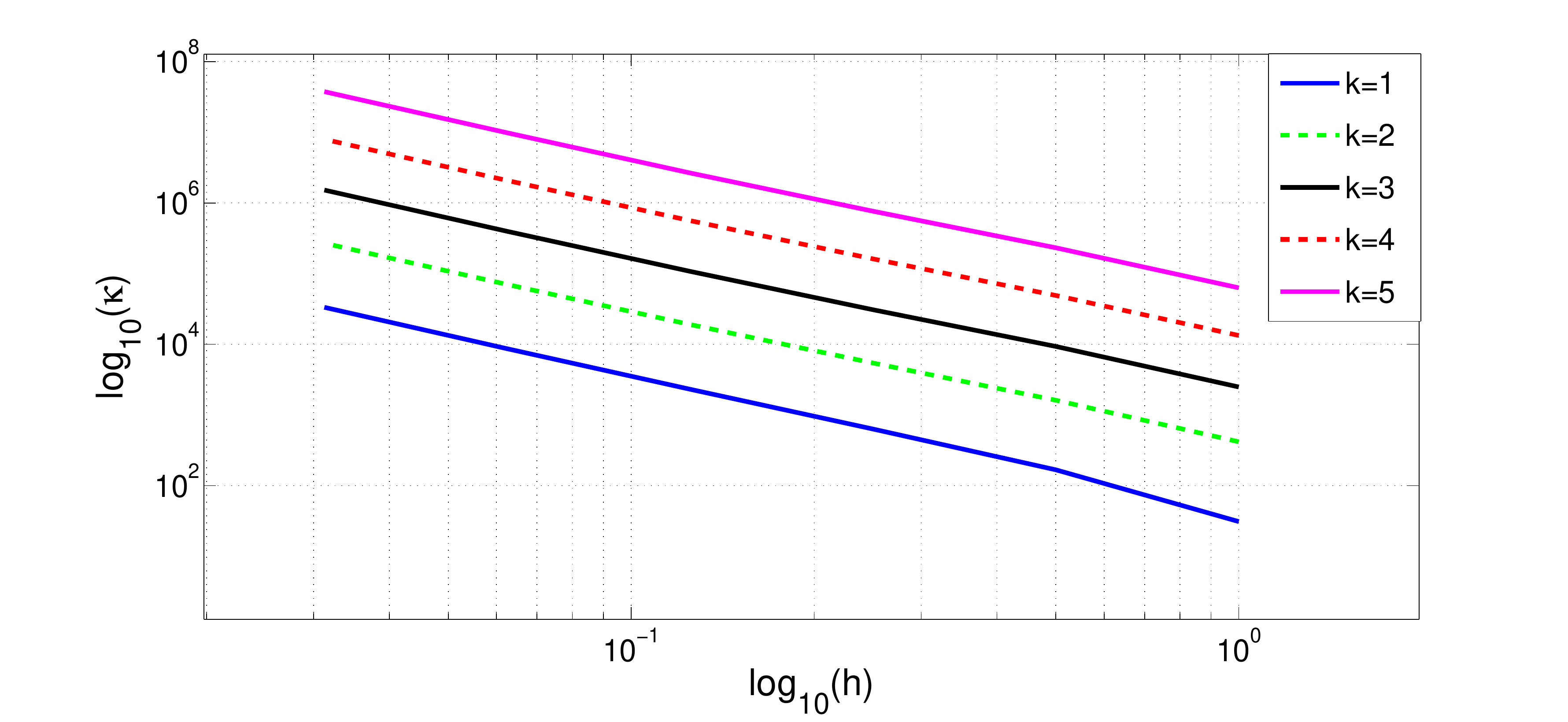}
\caption{Condition number of the stiffness matrix of the SIPG method as a function of mesh-size $h$ with different polynomial degree $k$ for the linear model problem (\ref{1}) with $r(u)=0$ }
\label{hmax_cond}
\end{figure}

Similarly, the condition number of the stiffness matrix increases with decreasing mesh-size and increasing order of the DG discretization for the linear diffusion-convection-reaction problem (3) with $r(u)=0$, (see Fig.\ref{hmax_cond}), which affects the efficiency of an iterative solver. Similar results can be found  in \citep{castillo12pdg} for the Poisson problem.

Besides the choice of a suitable basis polynomials, in this section, we describe also an efficient solution technique  for the large ill-conditioned linear systems arising from the linearization of the DG discretization. This technique is based on reordering of matrix elements and preconditioning.

\subsection{Matrix reordering \& block LU factorization }
Because the stiffness matrices obtained by DGFEM become ill-conditioned and dense for higher order DG elements \citep{ayuso09dgm}, several preconditioners are developed for the  efficient and accurate solution of linear diffusion-convection equations   \citep{antonetti02uad,georgoulis08npd}. We apply the matrix reordering and partitioning technique in
\citep{tari13rsp}, which uses the largest eigenvalue and corresponding eigenvector of the Laplacian matrix. This reordering allows us to obtain a partitioning and a preconditioner based on this partitioning. Since our matrices are non-symmetric, as the first step, we compute the symmetric structure by adding its transpose to itself. A symmetric, square and sparse matrix  could be represented as a graph where same index rows and columns are mapped into vertices  and  nonzeros of the sparse matrix are mapped into the edges of the graph. Since the matrix is symmetric, the corresponding graph is undirected. The  Laplacian matrix ($\mathcal{L}$) is, then, defined as follows
\[
\mathcal{L}(i,j)=
\begin{cases}
deg(v_i) &\text{if } i=j, \\
-1 &\text{if } i\neq j
\end{cases}
\]
in which the $deg(v_i)$ is the degree of the vertex i.  In this paper, the reordering we use is based on the unweighted Laplacian matrix given above. If the graph contains only one connected component, the eigenvalues of  $\mathcal{L}$ are
$0 = \lambda_1 < \lambda_2 \leq \lambda_3 \leq ... \leq \lambda_n$, otherwise there are as many zero eigenvalues as the number of connected components.

Certain eigenvalues and corresponding eigenvectors  of the Laplacian matrix have been studied extensively. Most notably the second nontrivial  eigenvalue of the Laplacian and the corresponding eigenvector known as the algebraic connectivity and the Fiedler vector of the graph \citep{fiedler73acg}.  Nodal domain theorem in \citep{fiedler75pen} shows that the eigenvectors corresponding to the eigenvalues other  than the first and the second smallest eigenvalue give us the connected components of the graph.  In \citep{barnard95sae},  the Fiedler vector for permuting the matrices to reduce the bandwidth is proposed. Reordering to obtain effective and scalable parallel banded preconditioners is proposed in \citep{manguoglu10wmo}. We use a sparse matrix reordering for partitioning and solving linear systems using the largest eigenvalue and the corresponding eigenvector of the Laplacian matrix.  Using this reordering, we show that one can reveal underlying structure of a sparse matrix. A simple Matlab code to find the reordered matrix and the permutation matrix can be found at  (\url{http://www.ceng.metu.edu.tr/~manguoglu/MatrixReorder.m})

To solve the discrete problem (\ref{ds}), we use the Newton-Raphson method. We start with a non-zero initial vector ${\bf U}^{(0)}$. The linear system arising from $i^{th}$-Newton-Raphson iteration step has the form $Jw^{(i)}=-R^{(i)}$, where $J$ is the Jacobian matrix to $R({\bf U}^{(0)})$ (i.e. $J=S+h'({\bf U}^{(0)})$ and it remains unchanged among the iteration steps), $w^{(i)}={\bf U}^{(i+1)}-{\bf U}^{(i)}$ is the Newton correction, and $R^{(i)}$ denotes the residual of the system at ${\bf U}^i$ ($R^{(i)}=R({\bf U}^{(i)})$). We construct a permutation matrix $P$ using the matrix reordering technique described above, applied to the sparse matrix J. Then, we apply the permutation matrix $P$ to obtain the permuted system $Nw=d$ where $N=PJP^T$, $w=Pw^{(i)}$  and  $d=-PR^{(i)}$. After solving the permuted system, the solution of the unpermuted linear system can be obtained by applying the inverse permutation, $w^{(i)} = P^T w$. Given a sparse linear system of equations $Nw=d$, after reordering, one way to solve this system is via block LU factorization. Suppose, the permuted matrix $N$, the right hand side $d$ and the solution $w$  is partitioned as follows:
\begin{equation*}
\begin{pmatrix}
A& B   \\
C^T     & D  \\
\end{pmatrix}
\begin{pmatrix}
w_1  \\
w_2  \\
\end{pmatrix}
=
\begin{pmatrix}
d_1  \\
d_2  \\
\end{pmatrix}
\end{equation*}
A block LU factorization of the coefficient matrix can be given as
\begin{equation*}
\begin{pmatrix}
A&  B   \\
C^T     & D  \\
\end{pmatrix}
=
\begin{pmatrix}
A&  0   \\
C^T     & S  \\
\end{pmatrix}
\begin{pmatrix}
I &  U  \\
0     & I  \\
\end{pmatrix}
\end{equation*}
where $U=A^{-1}B $ and  $S=D - C^TA^{-1}B$, also known as the Schur complement matrix. If the cost can be amortized, one can form $U$ and $S$ once and use them for solving linear systems with the same coefficient matrix and different right hand sides.  After this factorization, there are various approaches that one can take to solve the  system. One way is to solve the system via block backward and forward substitution, by first solving the linear system $A t=d_1$, and then solving the Schur complement system $Sw_2 = d_2 -C^T t$ and obtaining $w_1=t-Uw_2$. This method is summarized in Algorithm \ref{algorithm:linear}.

\begin{algorithm}
\caption{Algorithm for solving the linear system after reordering}
\textbf{Input:} The coefficient matrix: $\begin{pmatrix}
A&  B   \\
C^T     & D  \\
\end{pmatrix} $ and the right hand side:
$\begin{pmatrix}
d_1  \\
d_2  \\
\end{pmatrix}$
 \\
\textbf{Output:} The solution vector:
$\begin{pmatrix}
w_1  \\
w_2  \\
\end{pmatrix}$
 \\
\begin{algorithmic}[1]
\STATE  solve  $A t=d_1$
\STATE solve $Sw_2 = d_2 -C^T t$
\STATE compute $w_1=t-Uw_2$
\end{algorithmic}
\label{algorithm:linear}
\end{algorithm}

 We note that this approach involves solving two linear systems of equations  with the coefficient matrix $A$ and $S$. These linear systems can be solved directly or iteratively which requires effective preconditioners. Other approaches could involve solving the system $Nw=d$ iteratively where the preconditioner could take many forms. There are many other techniques for solving block partitioned and saddle point linear systems, we refer the reader to \citep{benzi05nss} for a more detailed survey of some of these methods.

\section{Numerical results}\label{numerical}
In this section, we give several numerical examples demonstrating the effectiveness and accuracy of the DGAFEM for convection dominated non-linear diffusion-convection-reaction equations.

\subsection{Example with polynomial type non-linearity}
\label{ex1}
Our first example is taken from \citep{bause10sfe} with Dirichlet boundary condition on $\Omega =(0,1)^2$ with $\epsilon =10^{-6}$, ${\bf b}=\frac{1}{\sqrt{5}}(1,2)^T$, $\alpha =1$ and $r(u)=u^2$. The source function $f$ and Dirichlet boundary condition are chosen so that $u(x_1,x_2)=\frac{1}{2}\left(1-\tanh \frac{2x_1-x_2-0.25}{\sqrt{5\epsilon}}\right)$ is the exact solution. The problem is characterized by  an internal layer of thickness $\mathcal{O}
(\sqrt{\epsilon}\mid\ln \epsilon\mid)$ around $2x_1 - x_2 = \frac{1}{4}$.

\begin{figure}[htb]
\centering
\subfloat[d.o.f. 70716]{\includegraphics[width=0.5\textwidth]{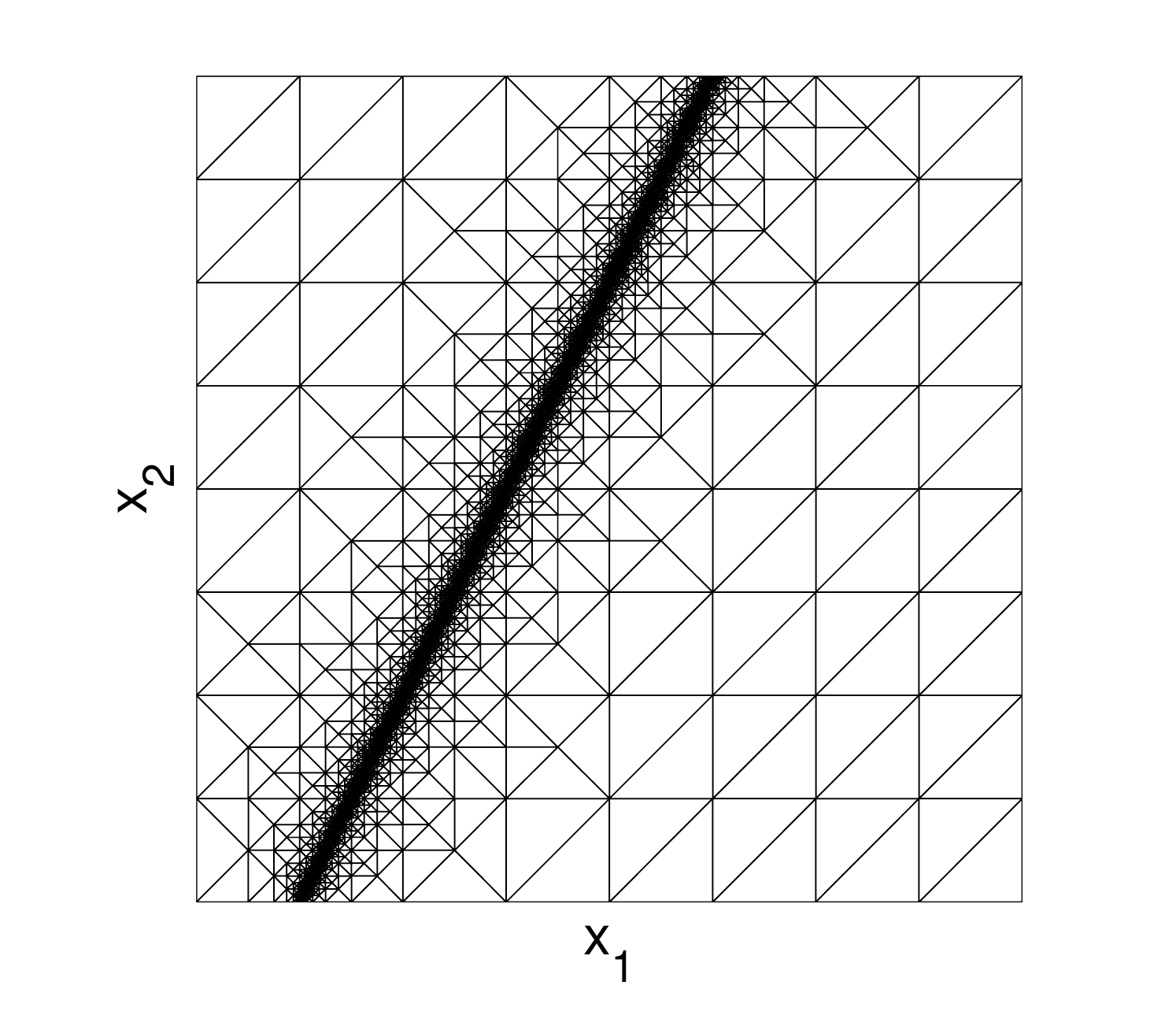}}
\caption{Example \ref{ex1}: Adaptive mesh. }
\label{ex1_mesh}
\end{figure}

\begin{figure}[htb]
\centering
\subfloat[d.o.f. 196608]{\includegraphics[width=0.6\textwidth]{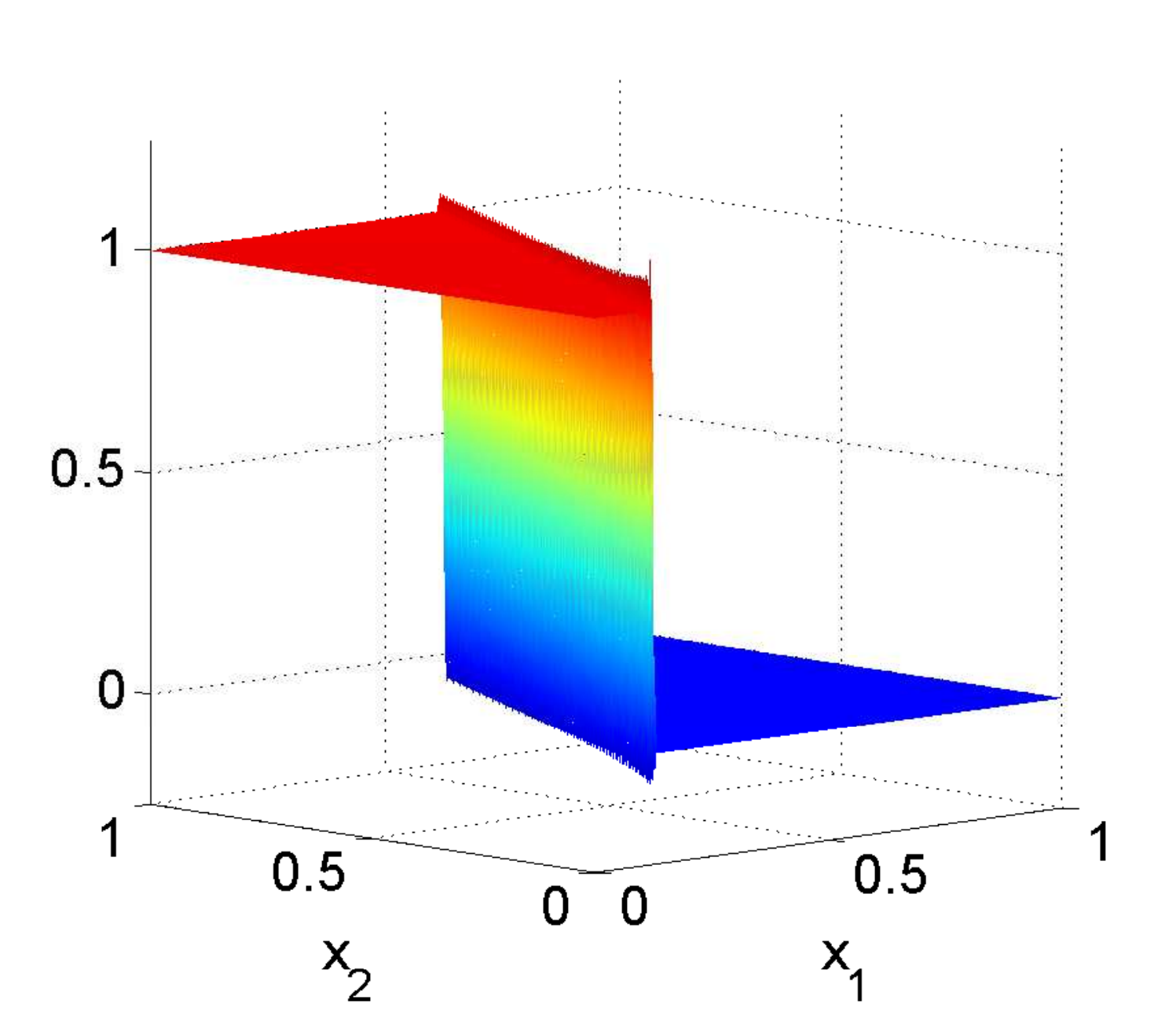}}
\subfloat[d.o.f. 70716]{\includegraphics[width=0.6\textwidth]{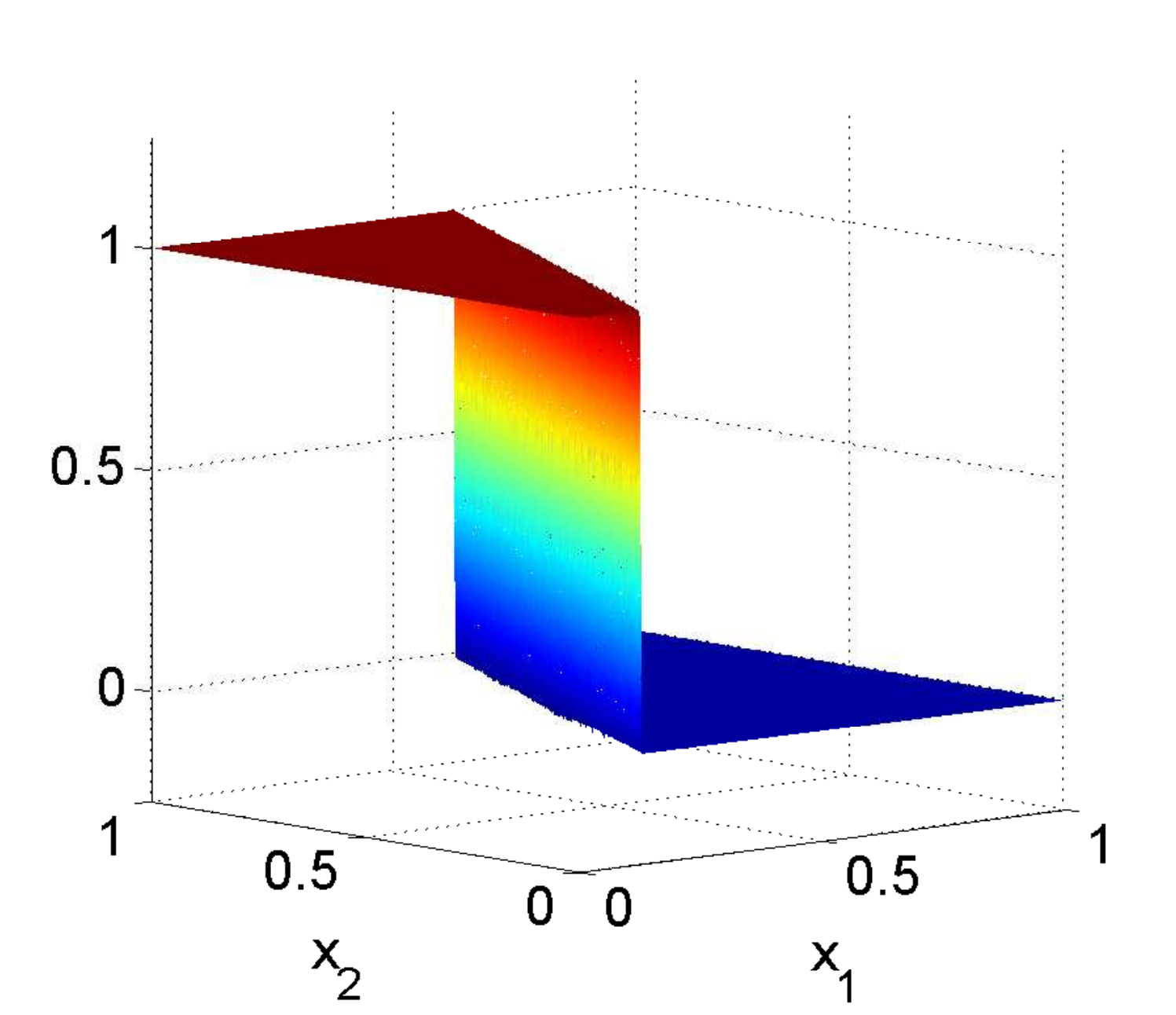}}
\caption{Example \ref{ex1}: Uniform(left) and adaptive(right) solutions, quadratic elements.}
\label{ex1_sol}
\end{figure}

\begin{figure}[htb]
\centering
\includegraphics[width=0.8\textwidth]{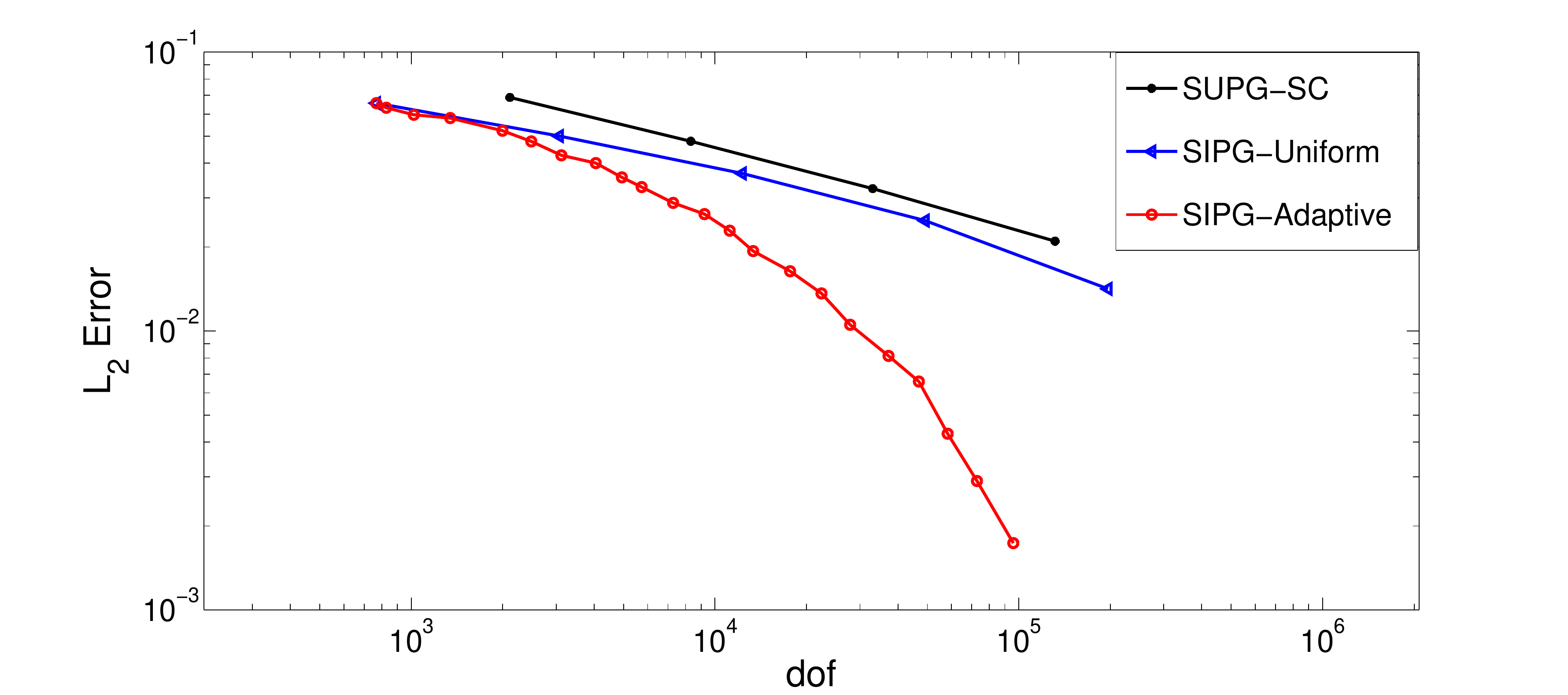}
\includegraphics[width=0.8\textwidth]{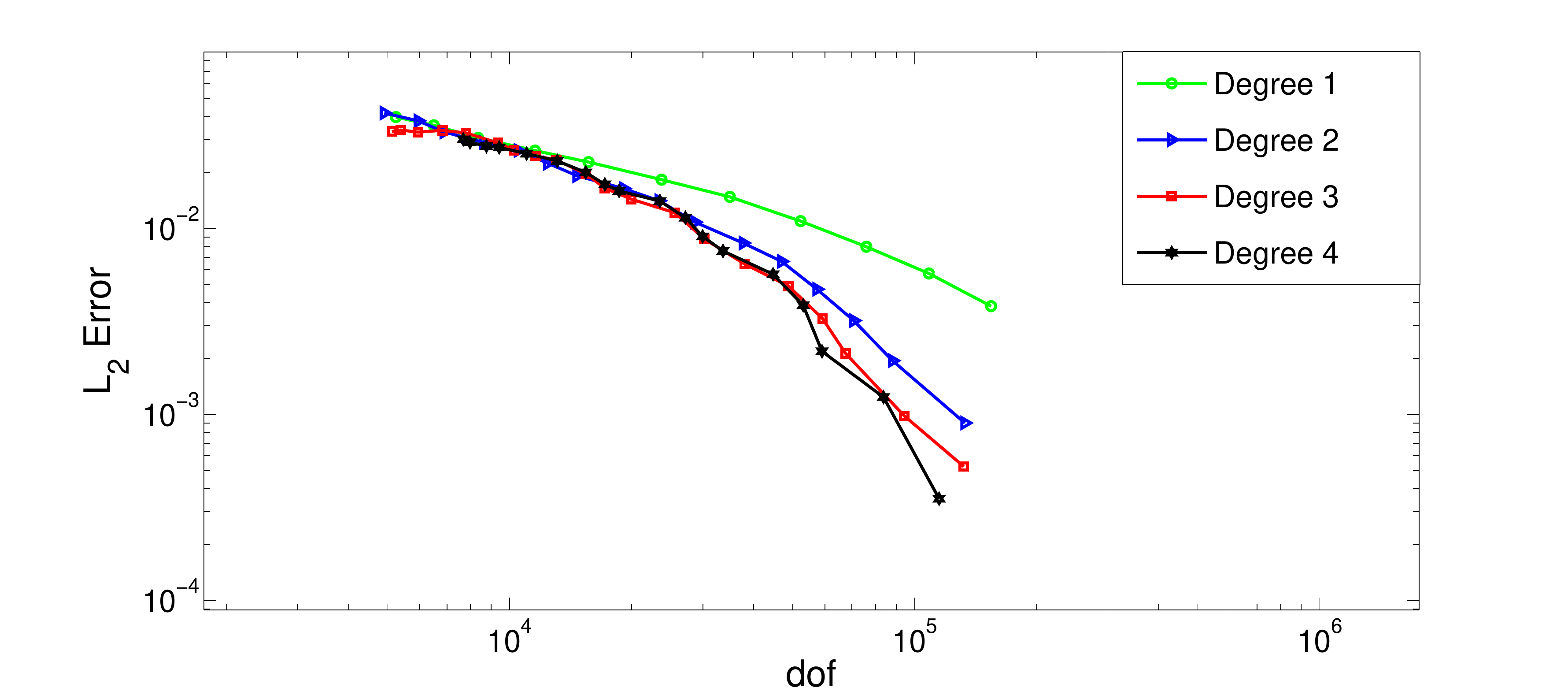}
\caption{Example \ref{ex1}: Global errors: comparison of the methods by quadratic elements(top), adaptive DG for polynomial degrees 1-4(bottom).}
\label{adapuniform}
\end{figure}

The mesh is locally refined by DGAFEM around the interior layer  (Fig.\ref{ex1_mesh}) and the spurious solutions are damped out in Fig.\ref{ex1_sol}, similar to the results as in  \citep{bause10sfe} using SUPG-SC, in \citep{yucel13dgf} with SIPG-SC. On adaptively and uniformly refined meshes, from the  Fig.\ref{adapuniform}, it can be clearly seen that the adaptive meshes reduce the substantial computing time. On uniform meshes, the SIPG is slightly more accurate as shown in \citep{yucel13dgf} than the the SUPG-SC in \citep{bause10sfe}. The error reduction by increasing degree of the polynomials is remarkable on finer adaptive meshes (Fig.\ref{adapuniform}, bottom).

In Table \ref{table_pol}, we give the results using the solution technique in Section \ref{linear} for the {\em BiCGStab}  with the stopping criterion as $ \| r_i\|_2/\| r_0\|_2\leq tol$ for $tol=10^{-7}$ ($r_i$ is the residual of the corresponding system at the $i^{th}$ iteration) applied to the unpermuted system and Schur complement system with and without preconditioning on the finest levels of uniformly ($4^{th}$ refinement level with dof 196608 and 32768 triangular elements) and adaptively ($17^{th}$ refinement level with dof 70716 and 11786 triangular elements) refined meshes. As a preconditioner, the incomplete LU factorization of the Schur complement matrix $S$ (ILU($S$)) is used for the linear system with the coefficient matrix $S$. The linear systems with the coefficient matrix $A$ are solved directly. Table \ref{table_pol} shows that solving the problem via the block LU factorization using the Schur complement system  with the preconditioner ILU($S$) is the fastest.

\begin{table}[htb]
\centering
\small\addtolength{\tabcolsep}{-5pt}
\begin{tabular}{l|l|l|l}
Linear Solver & \# Newton  & \# BiCGStab & Time \\
\hline
BiCGStab w/o prec. (Unpermuted)        & 10.8 (10.5) & 818 (757.5)   & 1389.3 (773.3) \\
BiCGStab w/ prec. $M_1$ (Permuted)     & 10.3 (10.3) & 1.5 (3)       & 423.1 (374.2) \\
BiCGStab w/ prec. $M_2$ (Permuted)     & 10.3 (10.3) & 1.5 (3)       & 416.8 (375.9) \\
Block LU + (BiCGStab w/o prec.)        & 10.3 (10.9) & 247.5 (315.5) & 270.9 (310.3)  \\
Block LU + (BiCGStab w/ prec. ilu(S) ) & 10.3 (10.9) & 19 (28.5)     & 140.9 (114.7)  \\
\hline
\end{tabular}
\caption{Example \ref{ex1}: Average number of Newton iterations, average number of {\em BiCGStab} iterations, total computation time in seconds corresponding to the uniformly refined (adaptively refined) mesh.}
\label{table_pol}
\end{table}

The time for  applying the permutation to obtain the reordered matrix and the permutation matrix $P$ takes $9.9$ seconds, whereas, it takes $0.13$ seconds to form the Schur complement matrix $S$ and $0.04$ seconds to compute the ILU($S$) on a PC with Intel Core-i7 processor and 8GB RAM
using the 64-bit version of Matlab-R2010a. We note that since the Jacobian matrix does not change during the non-linear
iterations, the permutation, the Schur complement matrix and ILU($S$) is computed only once for each run.

In all of the following results and figures, the Jacobian matrix $J$ is scaled  by a left Jacobi preconditioner before reordering to obtain a well conditioned matrix. The reordering procedure is applied to the scaled Jacobian matrix. Reordering times, which are included in the total computation time, for the uniform and adaptive refinements are $102.1$ seconds and $41.4$ seconds, respectively.

\begin{figure}[htb]
\centering
\includegraphics[width=0.8\textwidth]{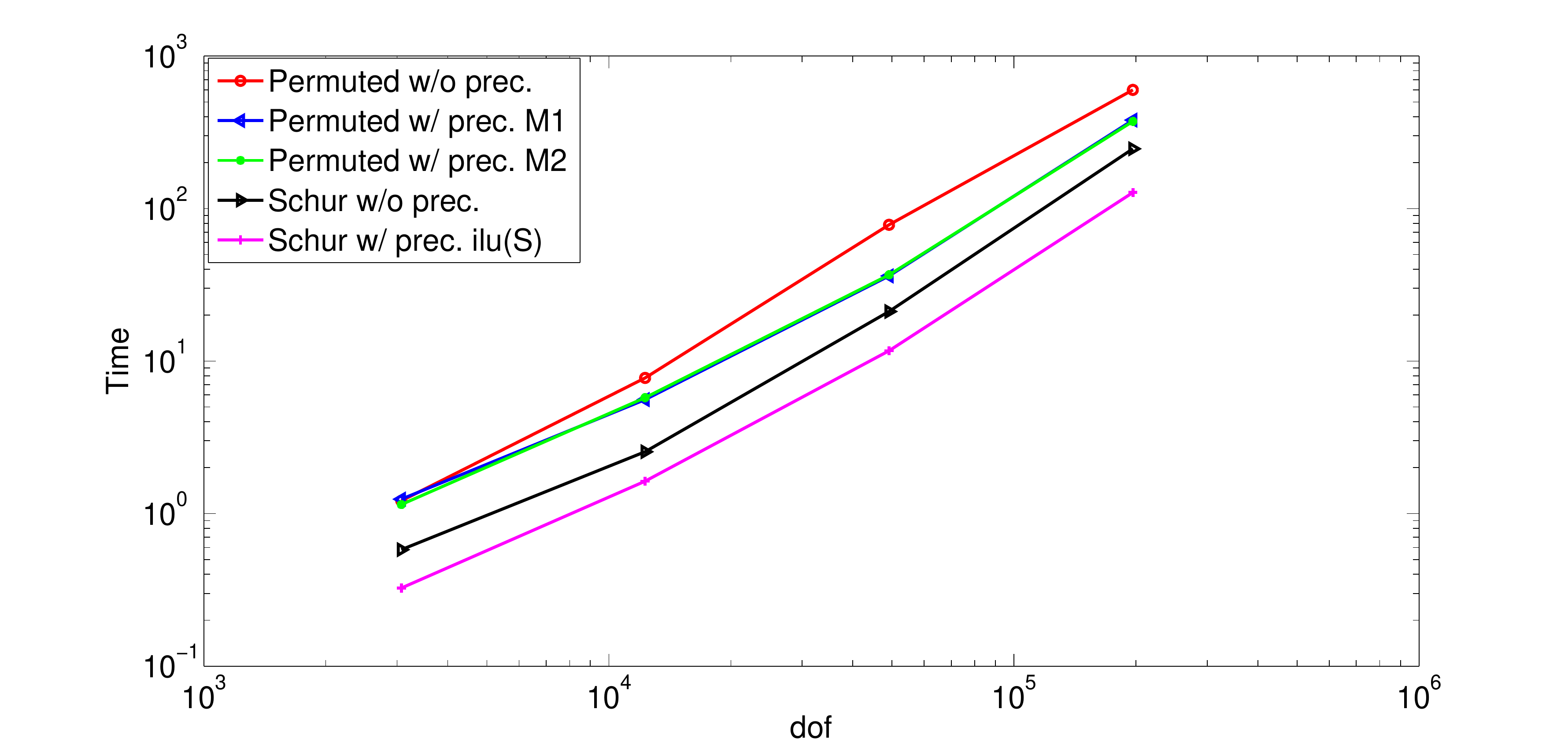}
\includegraphics[width=0.8\textwidth]{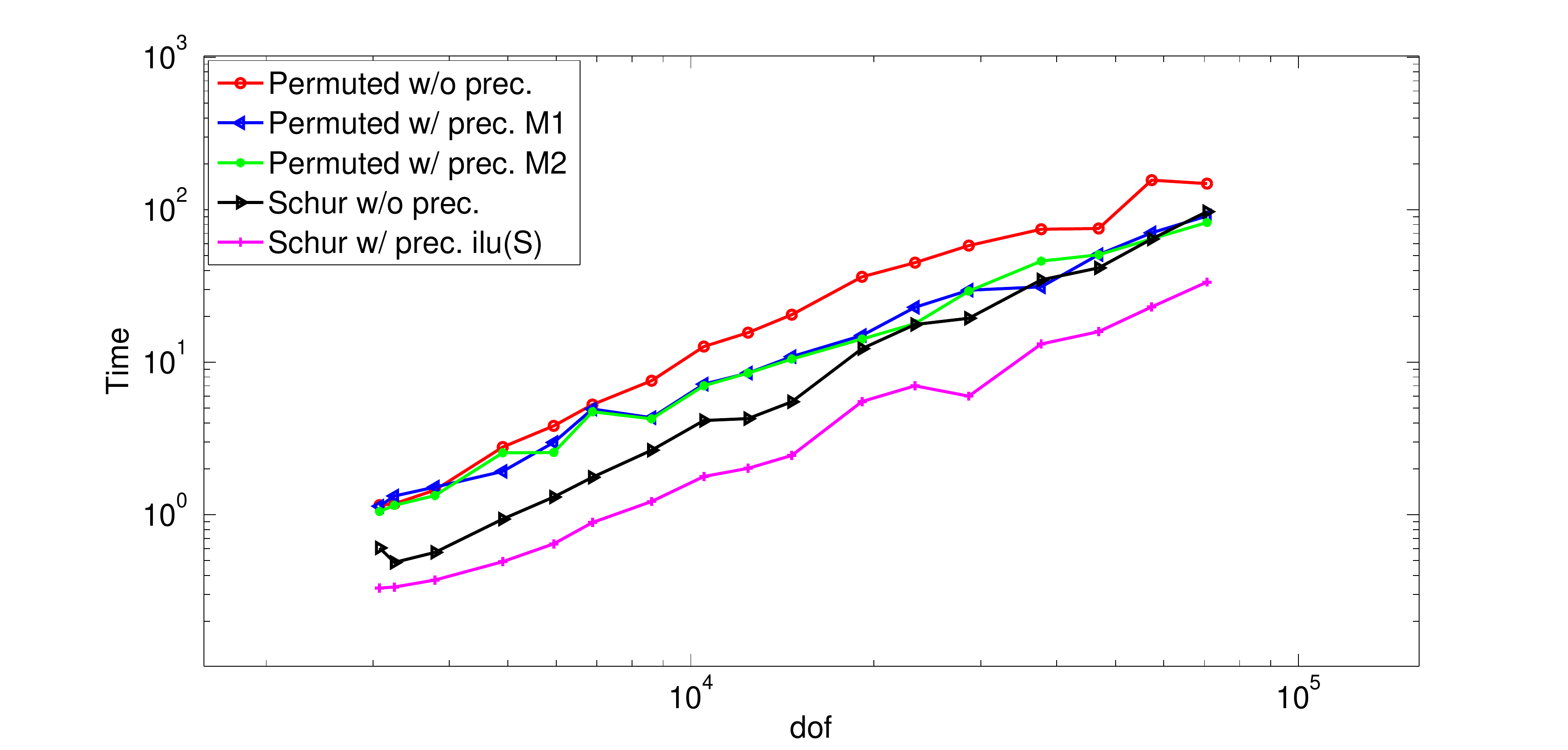}
\caption{Example \ref{ex1}: Computation time vs. dof: Uniform refinement (top) and adaptive refinement (bottom)}
\label{time_dof}
\end{figure}

\begin{figure}[htb]
\centering
\includegraphics[width=0.8\textwidth]{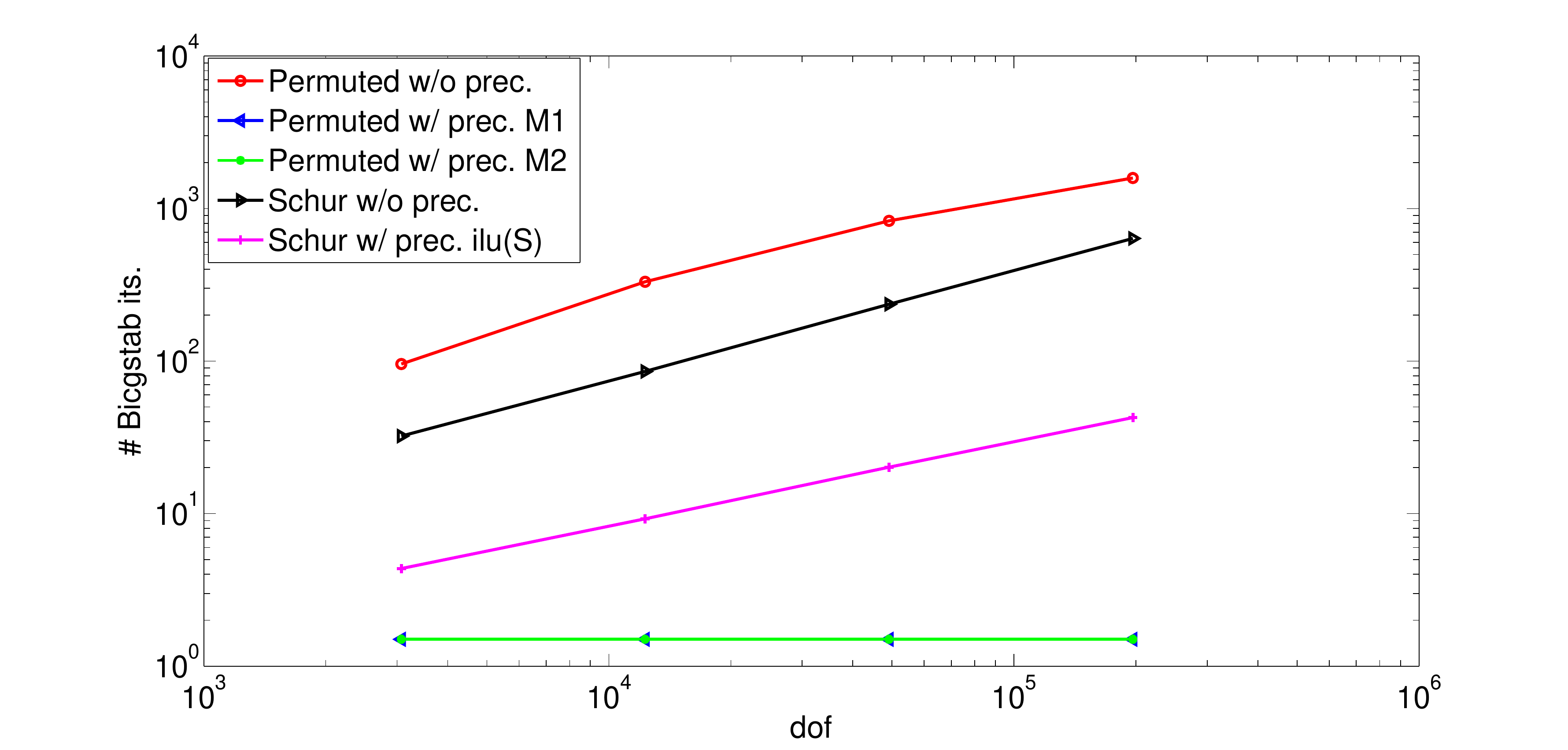}
\includegraphics[width=0.8\textwidth]{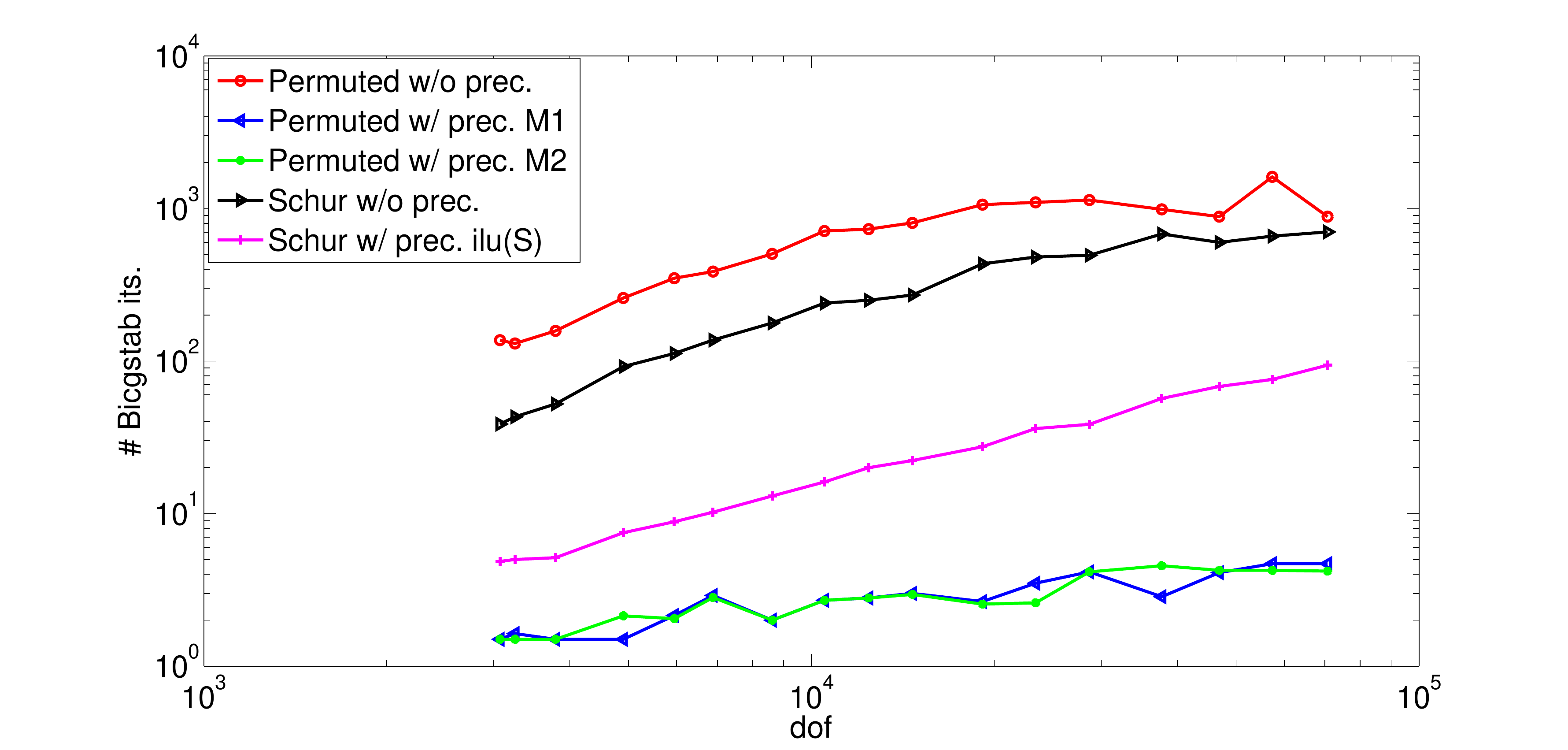}
\caption{Example \ref{ex1}: \# Average {\em BiCGStab} iterations vs. dof: Uniform refinement (top) and adaptive refinement (bottom)}
\label{its_dof}
\end{figure}

\begin{figure}[htb]
\centering
\includegraphics[width=0.8\textwidth]{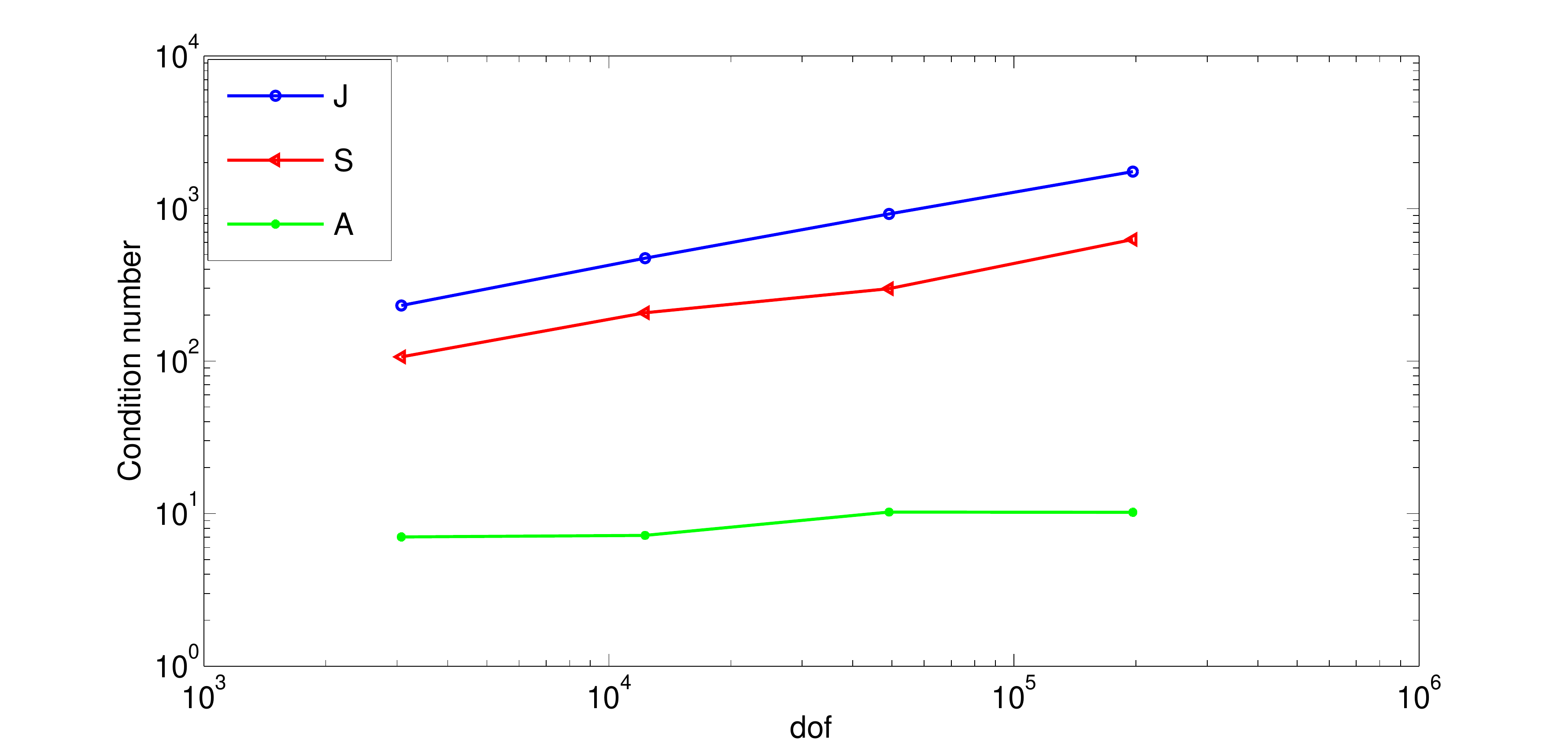}
\includegraphics[width=0.8\textwidth]{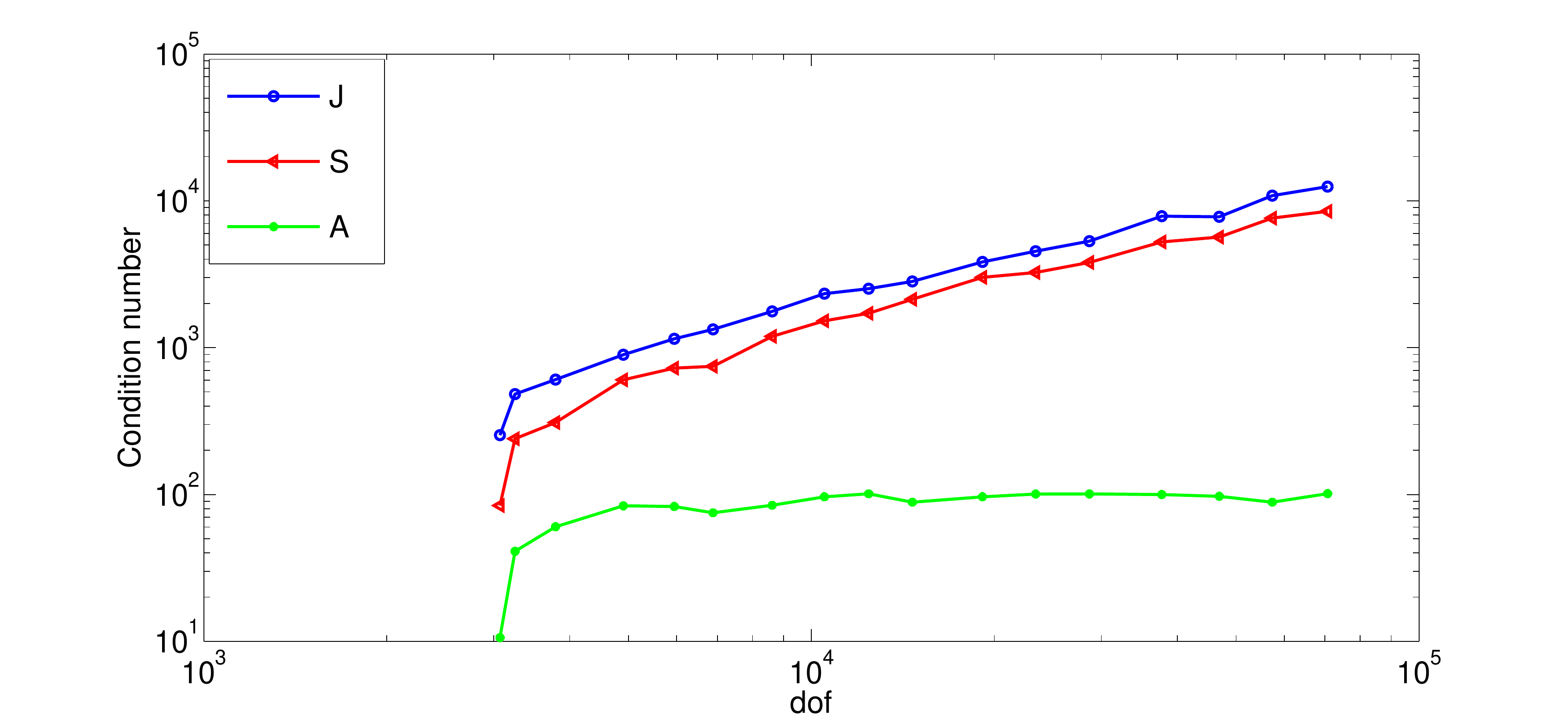}
\caption{Example \ref{ex1}: Condition number of the matrices J (unpermuted matrix), S (Schur complement matrix) and A (left top block of permuted matrix): Uniform refinement (top) and adaptive refinement (bottom).}
\label{ex1_cond}
\end{figure}

\begin{figure}[htb]
\centering
\includegraphics[width=0.9\textwidth]{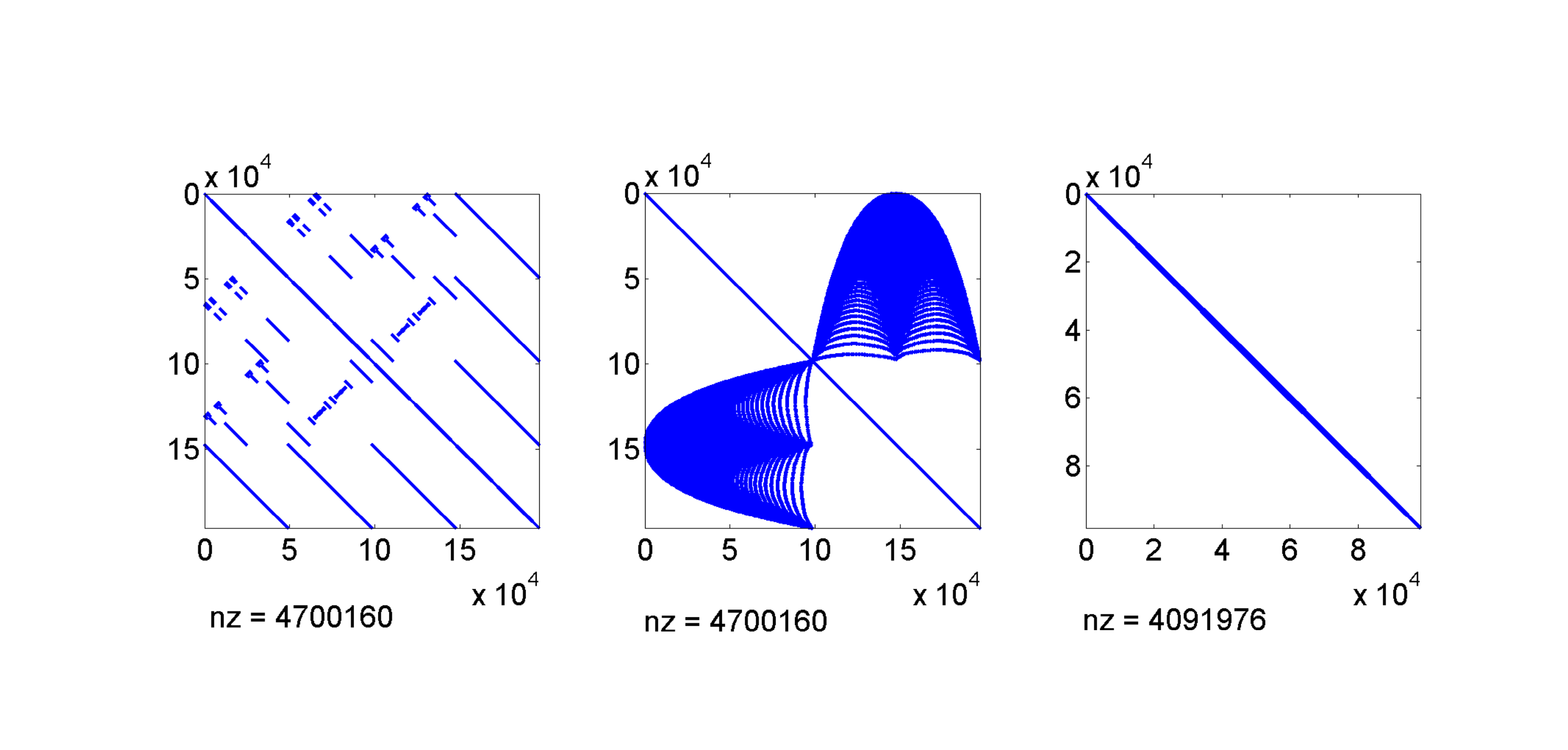}
\includegraphics[width=0.9\textwidth]{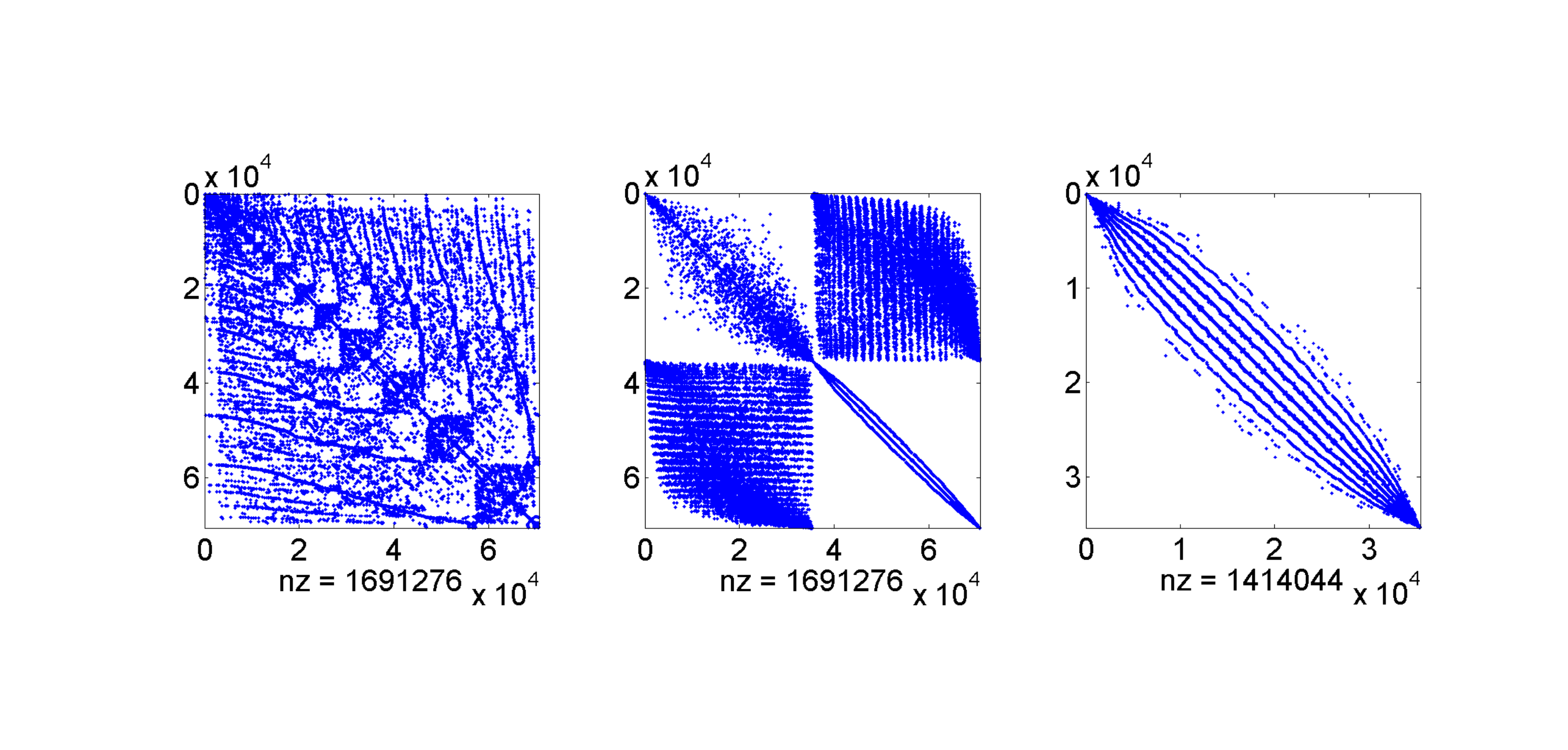}
\caption{Example \ref{ex1}: Sparsity patterns of the unpermuted (left), permuted (middle) and the Schur complement (right) matrices at the final refinement levels: Uniform refinement (top) with dof 196608 and adaptive refinement (bottom) with dof 70716.}
\label{ex1_spy}
\end{figure}

Fig.\ref{ex1_cond} shows the condition numbers of the Jacobian matrices $J$ of the original system, $S$ and $A$ of the block LU factorized system on the uniformly and adaptively refined meshes. The condition numbers of the coefficient matrix $A$ are almost constant for uniform refinement by different orders of DG discretizations and around $10$, whereas the condition number of  $S$ lower than of the Jacobian matrix $J$. This is due to the clustering of  nonzero elements around the diagonal (Fig.\ref{ex1_spy}) due to the matrix reordering. For adaptive refinement, Fig.\ref{ex1_cond}, bottom, we observe the same behavior, whereas the conditions numbers are larger of order one than for the uniform refinement. For comparison, we provide results by using BiCGStab with two block preconditioners. The preconditioning matrices $M_1$ and $M_2$ for the permuted full systems are given as
$$
M_1=
\begin{pmatrix}
A & 0\\
C^T & S
\end{pmatrix} \; , \qquad M_2=
\begin{pmatrix}
A & B\\
0 & S
\end{pmatrix}.
$$
Total number of iterations and time for different algorithms are given in Table \ref{table_pol}. Our proposed method where we compute the block LU factorization of the partitioned matrix and solve the system involving the Schur complement iteratively via preconditioned BiCGStab is the best in  terms of the total time compared to other methods for both uniform and adaptive refinement. In Fig.\ref{time_dof} and Fig.\ref{its_dof}, we present the total time and the average number of linear solver iterations, respectively, for uniform and adaptive refinements as the problem size has been increased. We observe that the proposed preconditioned linear solver has been the best in terms of time with a reasonable number of iterations for different problem sizes regardless refinement type.

\subsection{Example with Monod type non-linearity}
\label{ex2}
We consider the Monod type non-linearity in \citep{bause10sfe}:
$$
u-\nabla\cdot (\epsilon\nabla u)+\textbf{b}\cdot\nabla u-\frac{u}{1+u}=f
$$
on $\Omega =(0,1)^2$ with the convection field $\textbf{b}(x_1,x_2)=(-x_2,x_1)^T$, diffusion coefficient $\epsilon=10^{-6}$ and the source function $f=0$. The Dirichlet boundary condition is prescribed as $u(x_1,x_2)=1$ for $1/3\leq x_1 \leq 2/3, x_2=0$ and $u(x_1,x_2)=0$ on the remaining parts of the lower boundary as well as on the right and upper boundary. Moreover, $\frac{\partial u(x_1,x_2)}{\partial\textbf{n}}=0$ for $x_1=0, 0\leq x_2 \leq 1$ where \textbf{n} is the outer unit normal.

\begin{figure}[htb]
\centering
\subfloat{\includegraphics[width=0.5\textwidth]{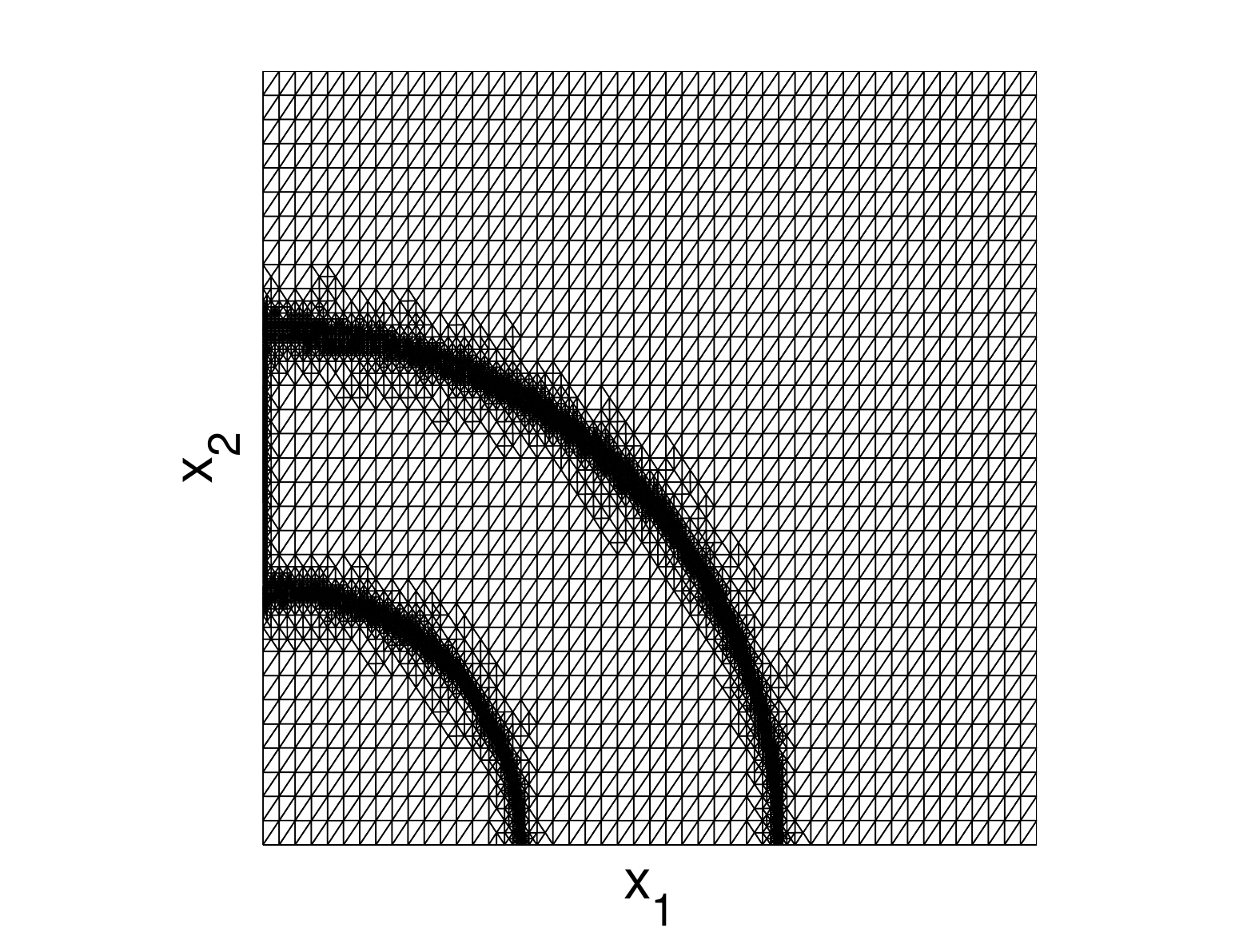}}
\subfloat{\includegraphics[width=0.65\textwidth]{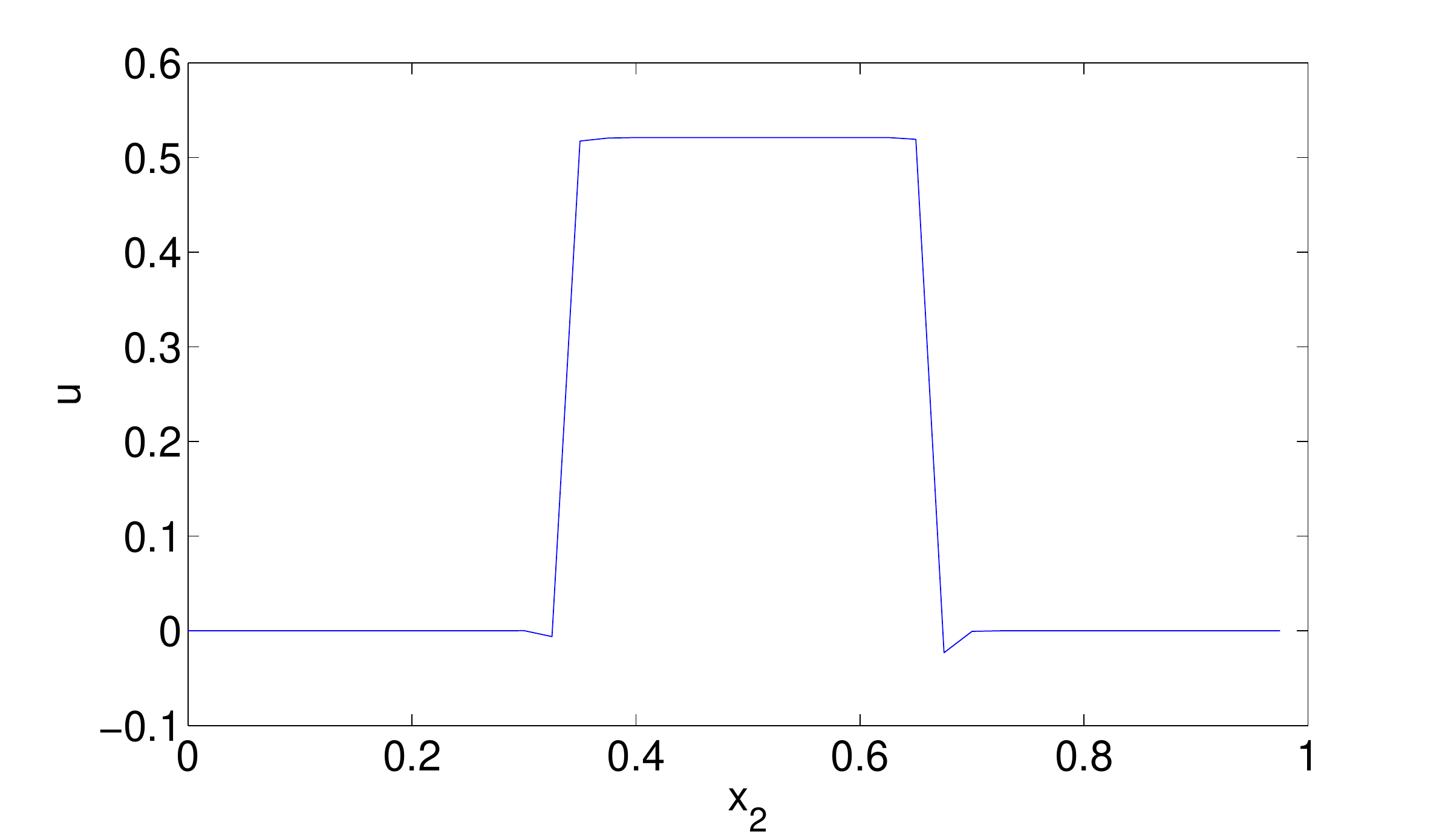}}

\centering
\subfloat{\includegraphics[width=0.6\textwidth]{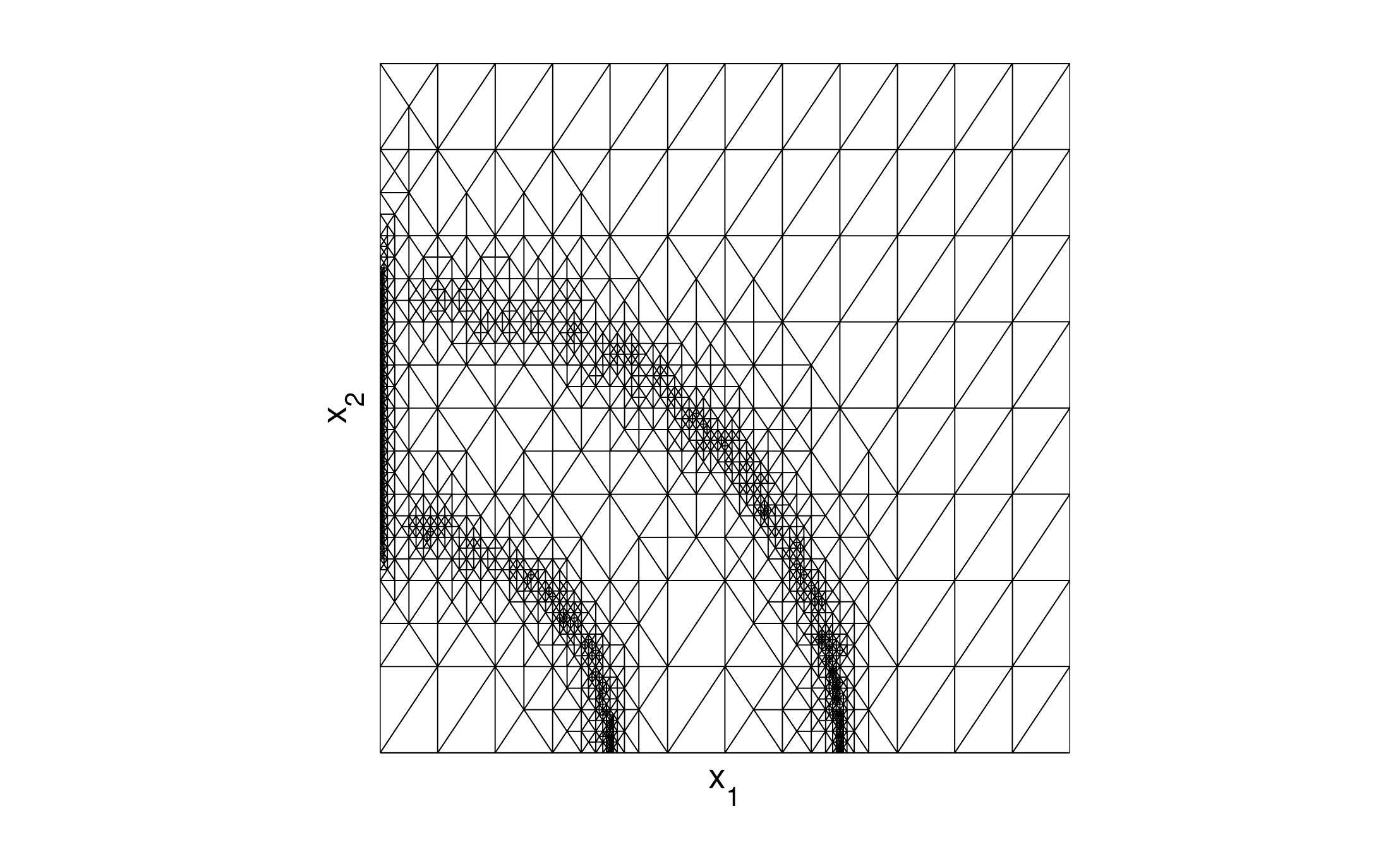}}
\subfloat{\includegraphics[width=0.7\textwidth]{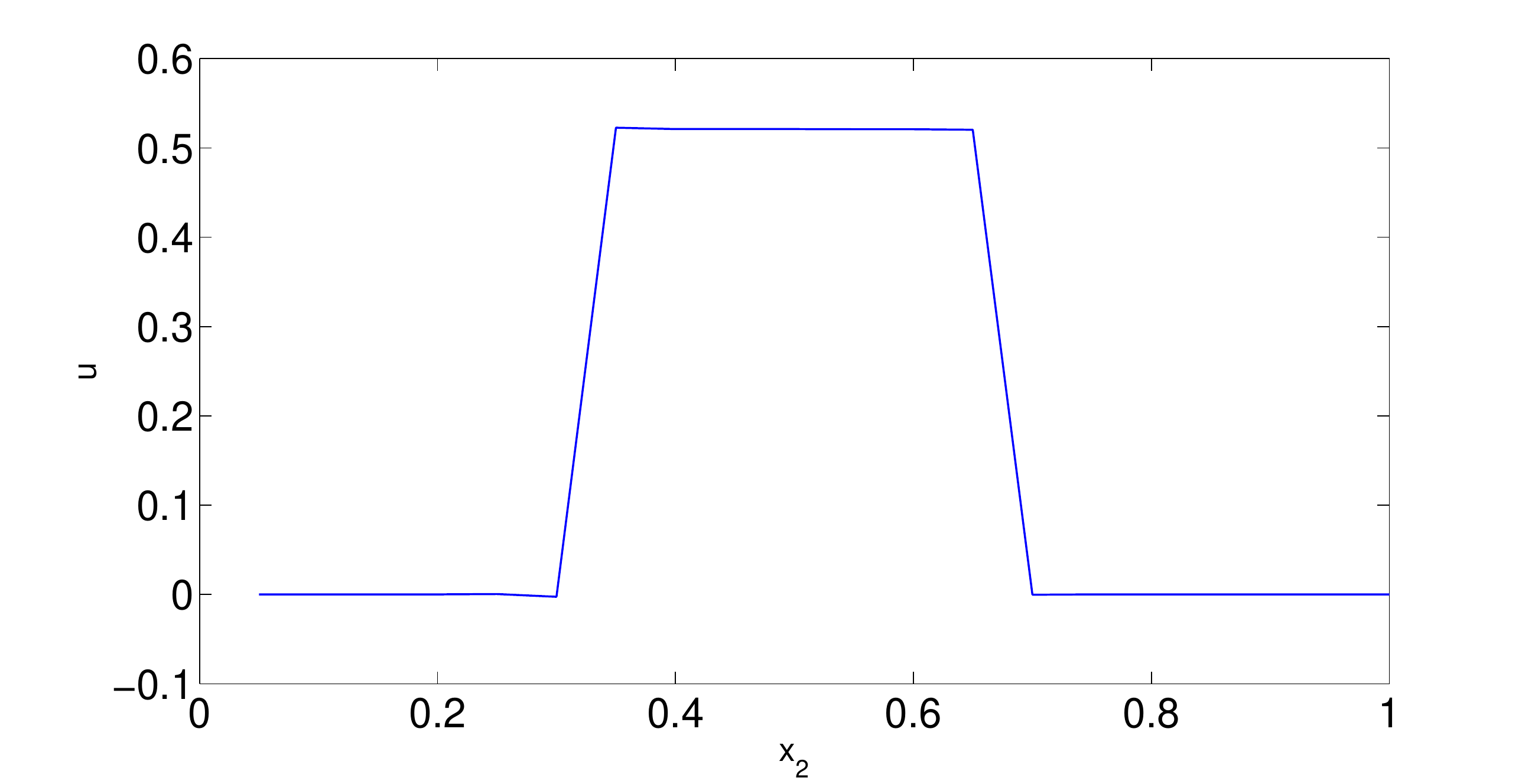}}
\caption{Example \ref{ex2}: Adaptive meshes (left) and the cross-section plots (right) of the solutions at the left outflow boundary by quadratic elements (top) with dof 82464 and quartic elements (bottom) with dof 80530 }.
\label{ex2_cross}
\end{figure}

There are both internal and boundary layers on the mesh (Fig.\ref{ex2_cross}, left), around them oscillations occur. Fig.\ref{ex2_cross}, right, shows that by DGAFEM, the oscillations are almost disappear, similar to the results in \citep{bause10sfe} for the  SUPG-SC and in \citep{yucel13dgf} for  SIPG-SC. Fig.\ref{ex2_cross}, left, shows that the adaptive process leads to correctly refined meshes. Moreover, by increasing polynomial degree of the basis functions ($k=4$), the oscillations are completely eliminated on the outflow boundary (Fig.\ref{ex2_cross}, bottom) and the sharp front is preserved. This is not the case for SUPG-SC \citep{bause10sfe}  and SIPG-SC \citep{yucel13dgf}, where still small oscillations are present.

As in case of polynomial non-linearity, Example 5.1, the block LU factorized system solved by BiCGStab with the preconditioner ILU(S) is the most efficient solver, with an average number of $7$ Newton iterations. The computing times  for the uniform refinement was 20.6 seconds, and 30.5 for the adaptive refinement.

\subsection{Example with Arrhenius type non-linearity}
\label{ex3}
Next example is the non-linear reaction for a two-component system in \citep{tezduyar86dcf}:
\begin{eqnarray*}
-\nabla\cdot (\epsilon\nabla u_1)+\textbf{b}\cdot\nabla u_1-100k_0u_2e^{\frac{-E}{Ru_1}}&=&0, \\
-\nabla\cdot (\epsilon\nabla u_2)+\textbf{b}\cdot\nabla u_2+k_0u_2e^{\frac{-E}{Ru_1}}&=&0
\end{eqnarray*}
on $\Omega =(0,1)^2$ with the convection field $\textbf{b}=(1-x_2^2,0)^T$, the diffusion constant $\epsilon=10^{-6}$, the reaction rate coefficient $k_0=3\times 10^8$ and the quotient of the activation energy to the gas constant $\frac{E}{R}=10^4$. The unknowns $u_1$ and $u_2$ represent the temperature of the system and the concentration of the reactant, respectively.

\begin{figure}[htb]
\centering
\subfloat{\includegraphics[width=0.5\textwidth]{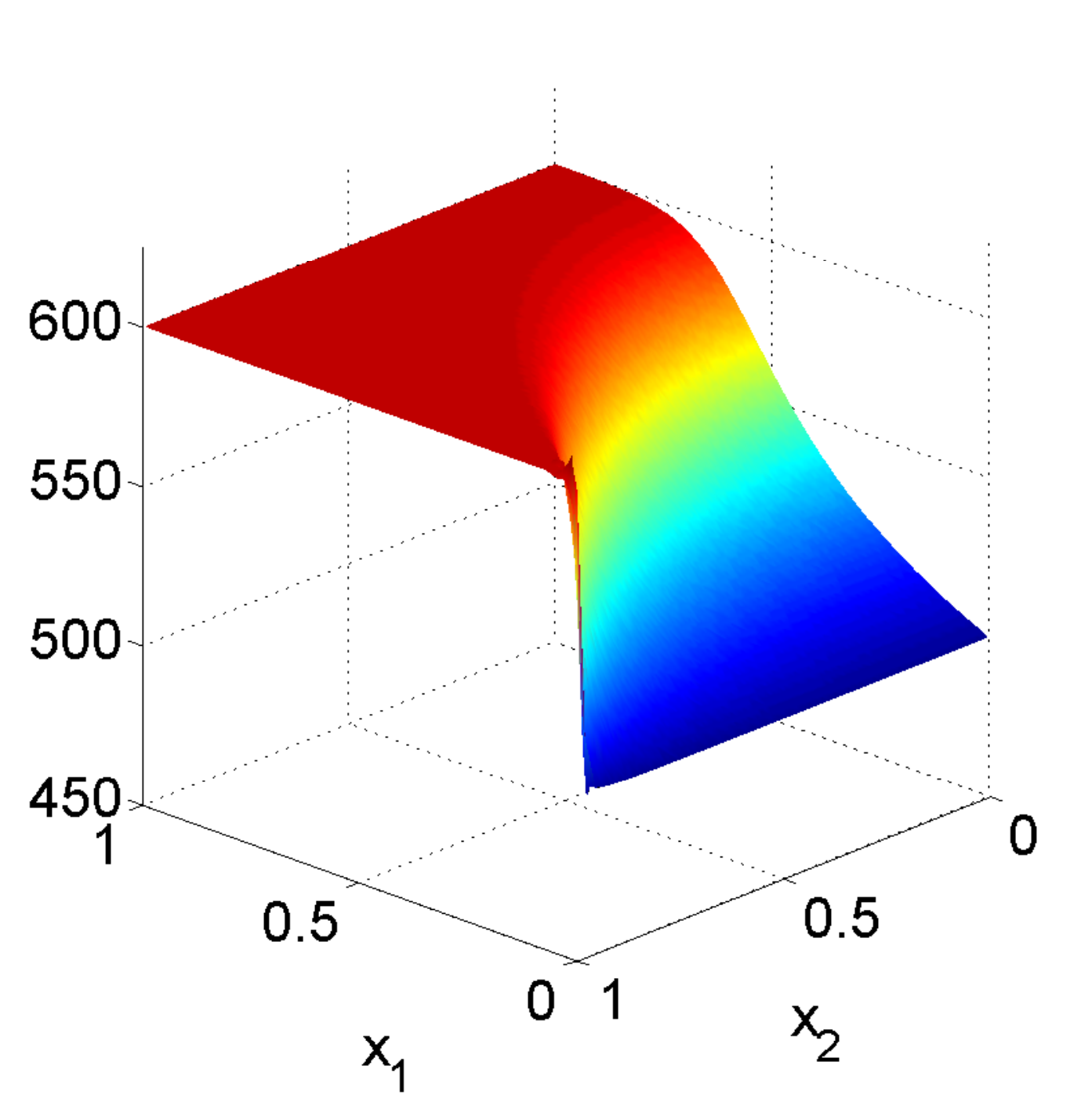}}
\subfloat{\includegraphics[width=0.5\textwidth]{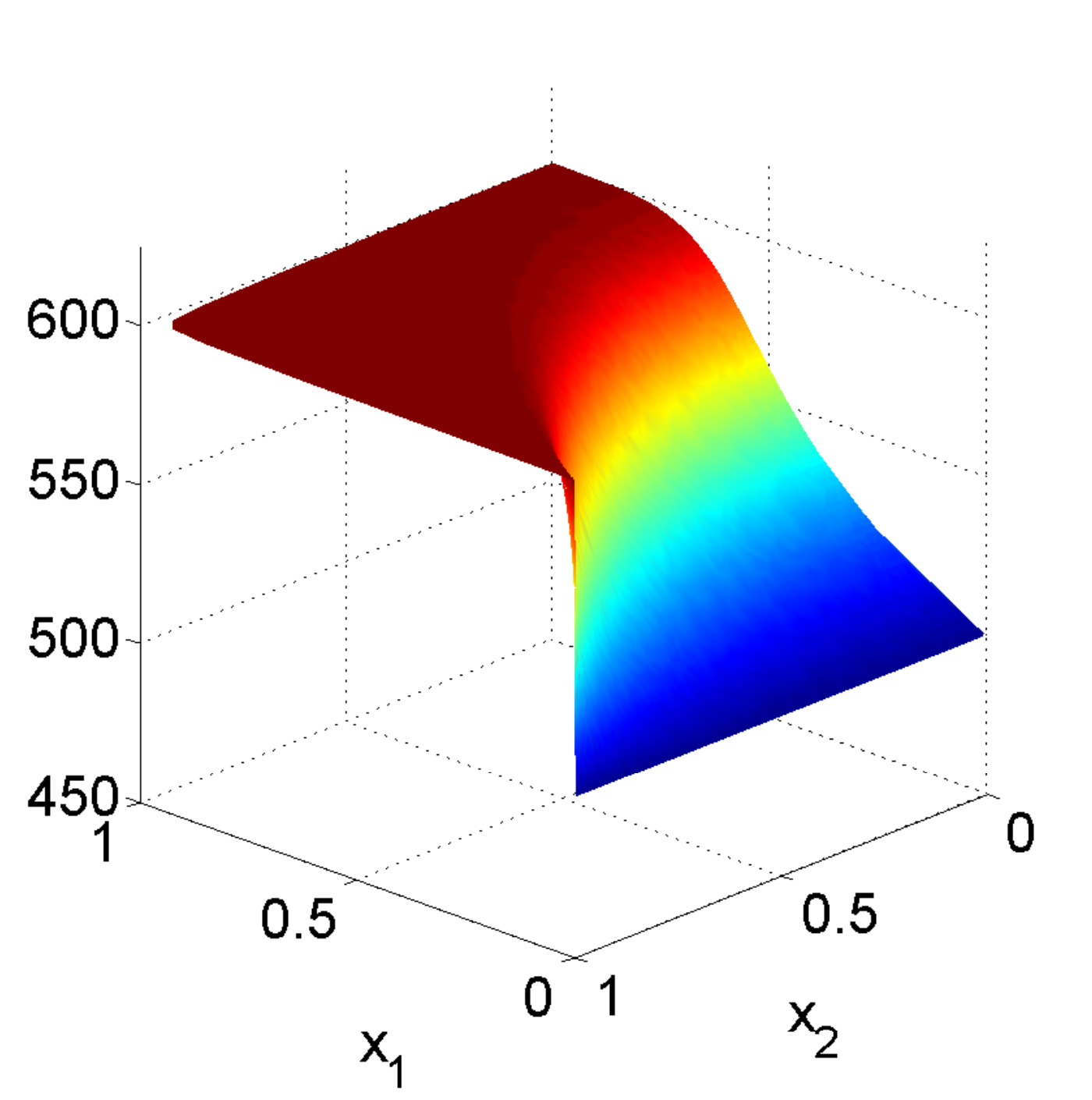}}

\centering
\subfloat{\includegraphics[width=0.5\textwidth]{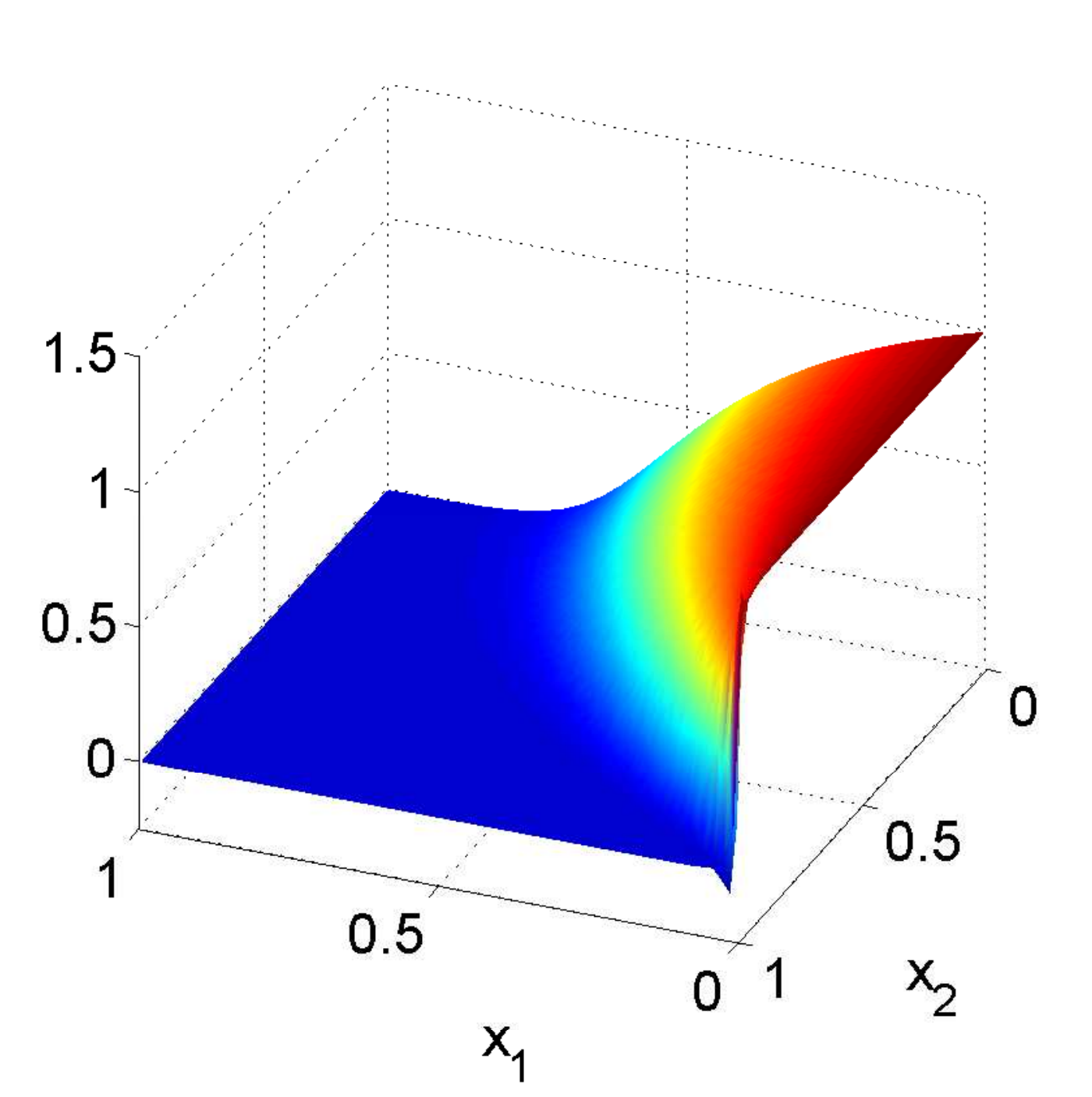}}
\subfloat{\includegraphics[width=0.5\textwidth]{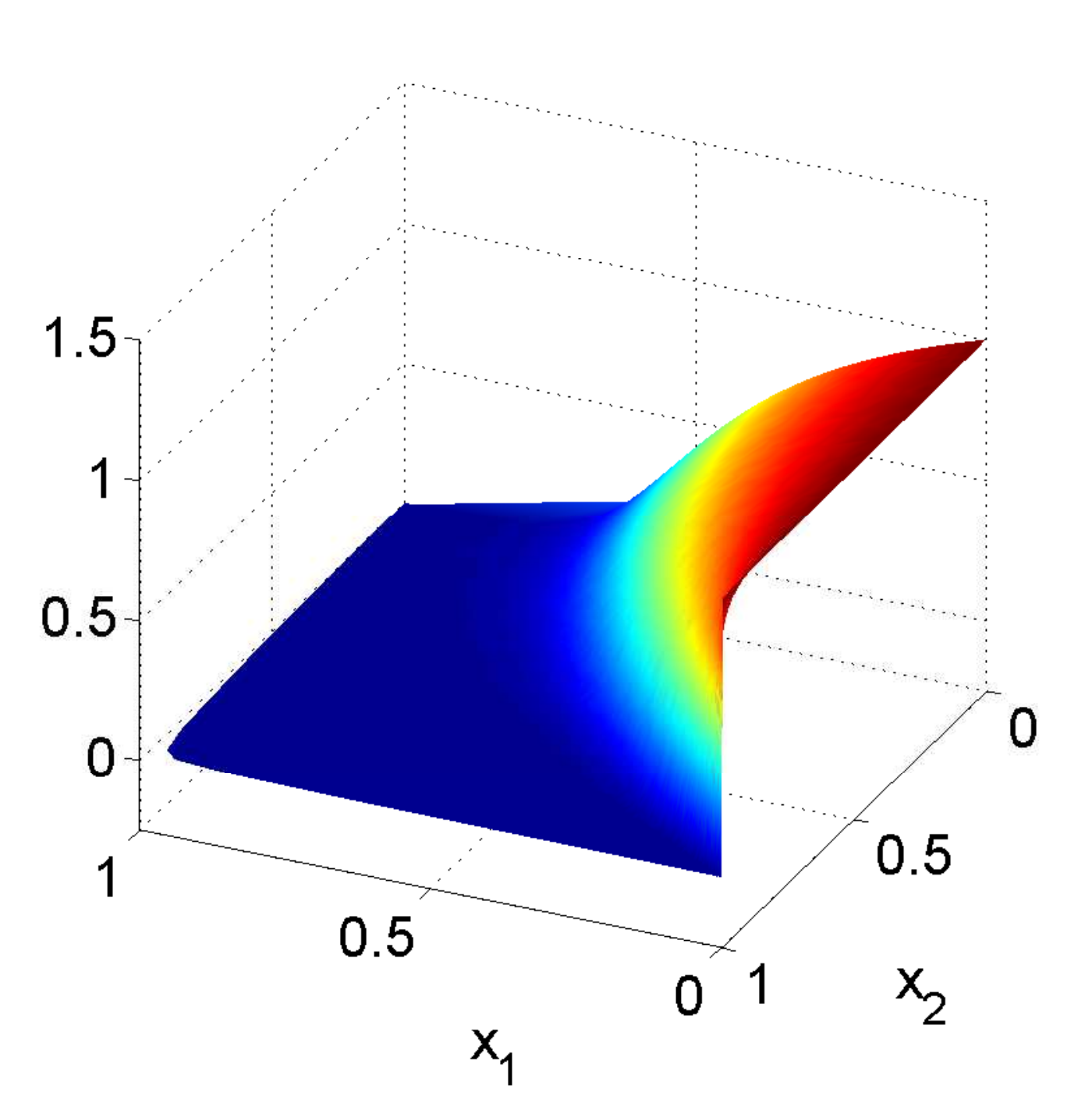}}
\caption{Example \ref{ex3}: Uniform(left) and adaptive(right) solutions to the temperature(top) and reactant(bottom), quadratic elements with dof 12288 for uniform refinement and with dof 6168 for adaptive refinement}.
\label{plot_arh}
\end{figure}

There are oscillations around the layers, even small, for the uniform refinement (Fig.\ref{plot_arh}, left)  as for SIPG-SC in \citep{yucel13dgf}. On the other hand, these oscillations are completely dumped out by DGAFEM with almost half of the dof used in the uniform refinement (Fig.\ref{plot_arh}, right).

The block LU factorization based algorithm with the preconditioner ILU(S) requires 10.5 seconds for the uniform and 24.4 seconds for the adaptive refinements. Matrix reordering and permutation took 2.44 seconds for the uniform  and  2.17 seconds for adaptive refinements, respectively.

\subsection{Coupled example with polynomial type non-linearity}
\label{ex4}
Our final problem is the modification of the non-stationary transport problem, Example 2, in \citep{bause13hof}. The problem is stated as the following:

\begin{eqnarray*}
\alpha u_1-\nabla\cdot (\epsilon\nabla u_1)+\textbf{b}\cdot\nabla u_1+50u_1^2u_2^2 &=& 0, \\
\alpha u_2-\nabla\cdot (\epsilon\nabla u_2)+\textbf{b}\cdot\nabla u_2++50u_1^2u_2^2 &=& 0
\end{eqnarray*}
on the rectangular domain $\Omega =(0,1)\times (0,2)$ with the convection field $\textbf{b}=(0,-1)^T$, the diffusion constant $\epsilon=10^{-10}$ and linear reaction constant $\alpha =0.1$. On the left, right and lower parts of the boundary of the domain, Neumann boundary conditions are prescribed. On the remaining part of the boundary, Dirichlet boundary conditions are chosen as
\[
u_1({\bf x})=
\begin{cases}
8(x_1-0.375) &\text{for } 0.375 < x_1 \leq 0.5 ,\\
-8(x_1-0.625) &\text{for } 0.5 < x_1 \leq 0.625 ,\\
0 &\text{otherwise}
\end{cases}
\]

\[
u_2({\bf x})=
\begin{cases}
8(x_1-0.125) &\text{for } 0.125 \leq x_1 \leq 0.25 ,\\
-8(x_1-0.375) &\text{for } 0.25 < x_1 \leq 0.375 ,\\
8(x_1-0.625) &\text{for } 0.625 \leq x_1 \leq 0.75 ,\\
-8(x_1-0.875) &\text{for } 0.75 < x_1 \leq 0.875 ,\\
0 &\text{otherwise}
\end{cases}
\]

\begin{figure}[htb]
\centering
\subfloat{\includegraphics[width=0.5\textwidth]{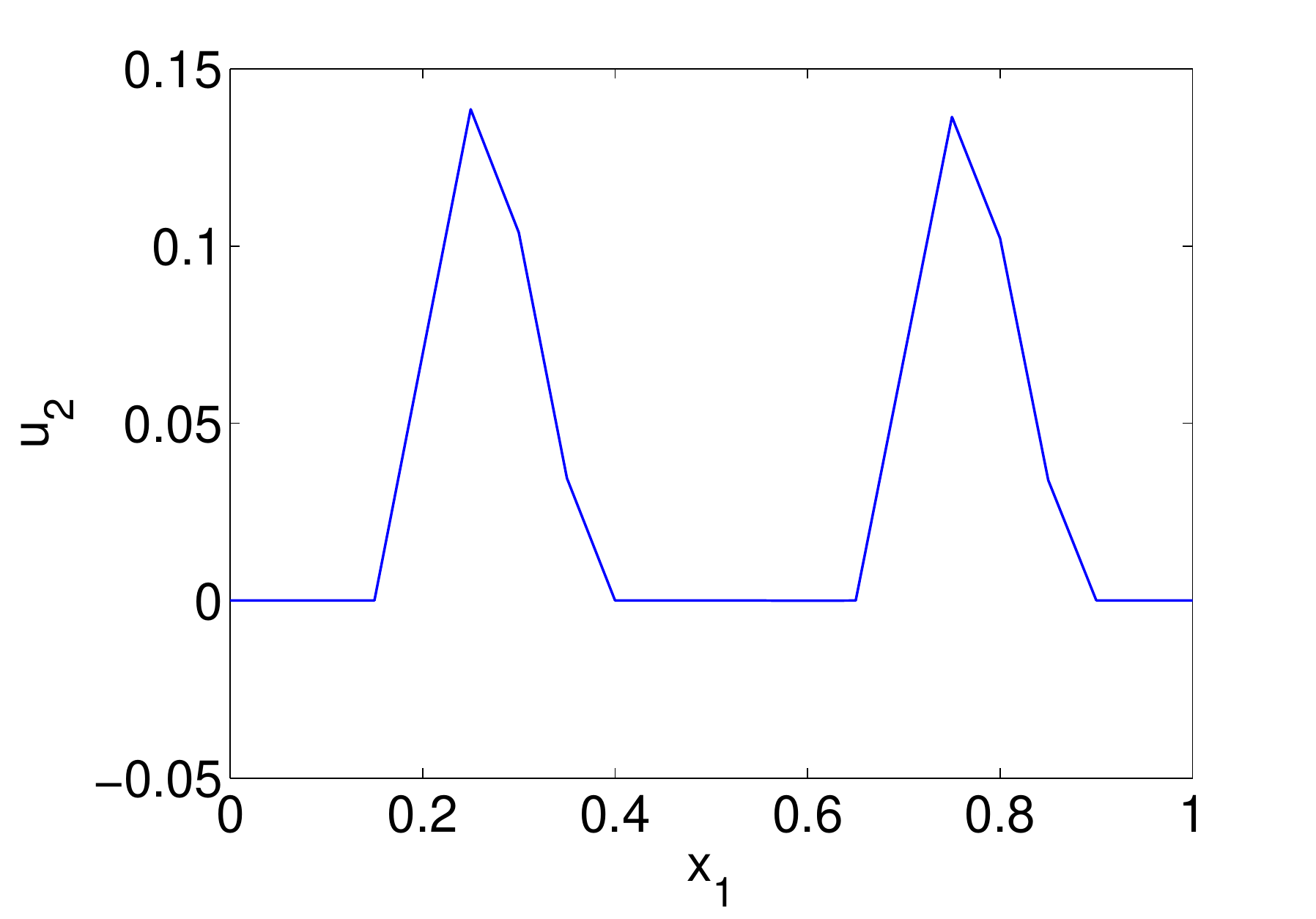}}
\subfloat{\includegraphics[width=0.5\textwidth]{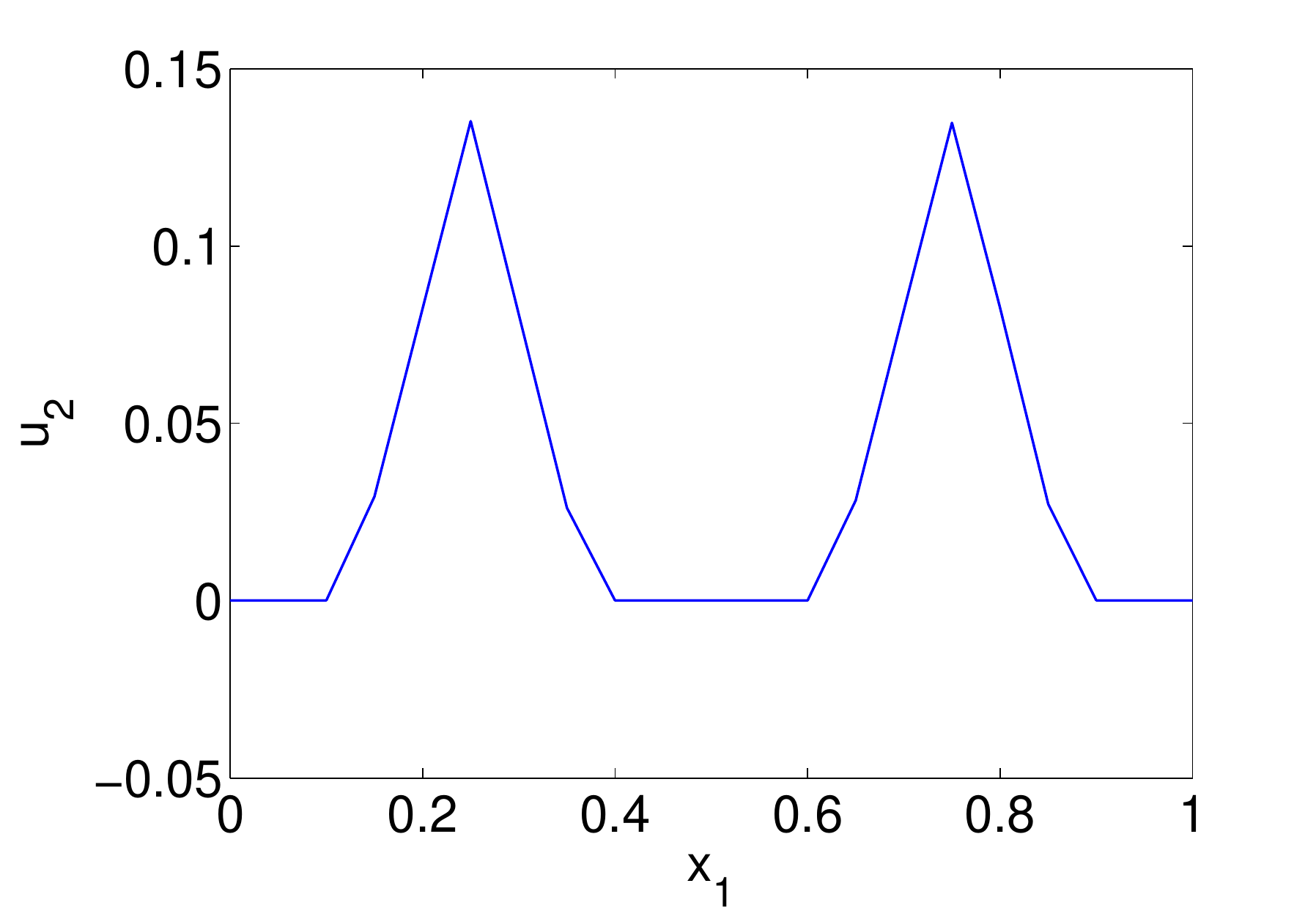}}
\caption{Example \ref{ex4}: Uniformly(left) and adaptively(right) obtained cross-section plots  on the outflow boundary for the component $u_2$, quartic elements with dof 61440 for uniform refinement and with dof 33690 for adaptive refinement.}
\label{ex4_plot}
\end{figure}

\begin{figure}[htb]
\centering
\includegraphics[width=0.5\textwidth]{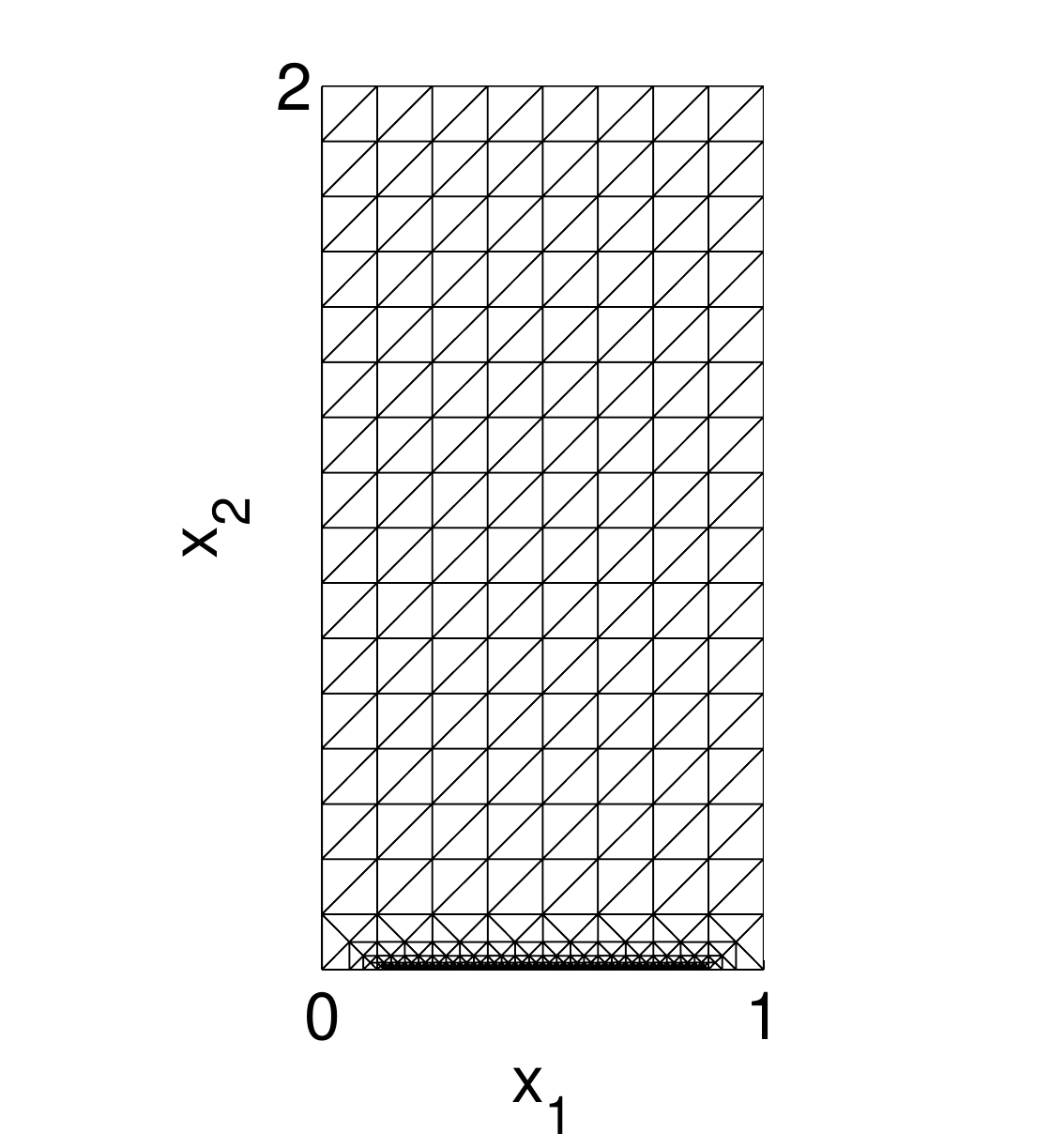}
\caption{Example \ref{ex4}: Adaptive mesh, quartic elements with dof 33690.}
\label{ex4_mesh}
\end{figure}

There is a boundary layer on the outflow boundary, Fig.\ref{ex4_mesh}. Fig.\ref{ex4_plot} shows that oscillations are almost damped using DGAFEM approximations, similar to those results in \citep{bause13hof} using SUPG-SC.  It can be seen from Fig.\ref{ex4_mesh} that the mesh is correctly refined by DGAFEM near the boundary layer.

\section{Conclusions}
We have shown that using DGAFEM with the sparse linear solver  is an efficient method for solving non-linear convection dominated problems accurately and avoids the design of the parameters in the shock capturing technique as for the SUPG-SC and DG-SC methods. The numerical examples demonstrate that DGAFEM allows to capture the interior and boundary layers very sharply without any significant oscillation. As a future study, we will apply space-time adaptive DG methods for time-dependent convection dominated non-linear diffusion-convection-reaction equations.

\section*{Acknowledgment}
The authors would like to thank the reviewer for the comments and suggestion that help improve the manuscript.
This work has been partially supported by Turkish Academy of Sciences
Distinguished Young Scientist Award TUBA-GEBIP/2012-19, T\"UBITAK Career
Award EEAG111E238 and METU BAP-07-05-2013-004.

\end{document}